\documentclass[reqno,a4paper]{amsart}
	
\usepackage[english]{babel}
\usepackage[utf8]{inputenc}

\usepackage{amsfonts,amsmath,amssymb,amsthm}
\usepackage[centering]{geometry}
\usepackage[hidelinks]{hyperref}

\usepackage{csquotes}
\usepackage[backend=biber,sorting=nyt]{biblatex}
\usepackage{mathtools}
\mathtoolsset{showonlyrefs=true}

\author[A. Pinoy]{Alan Pinoy}
\address{Institutionen för Matematik\\
Kungliga Tekniska Högskolan\\
100 44 Stockholm\\
Sweden}
\email{pinoy{@}kth{.}se}

\title[CR compactification for ALCH almost Hermitian manifolds]{CR compactification for asymptotically locally complex hyperbolic almost Hermitian manifolds}

\theoremstyle{plain}
\newtheorem{corollary}{Corollary}[section]
\newtheorem{lemma}[corollary]{Lemma}
\newtheorem{proposition}[corollary]{Proposition}

\newtheorem*{mainthm}{Main Theorem}
\newtheorem{thm}[corollary]{Theorem}
\newtheorem*{thm*}{Theorem}

\theoremstyle{definition}
\newtheorem{definition}[corollary]{Definition}

\theoremstyle{remark}
\newtheorem{remark}[corollary]{Remark}

\newcommand{\C}{\mathbb{C}}
\newcommand{\R}{\mathbb{R}}
\newcommand{\dK}{\partial \! K}

\DeclareMathOperator{\dmesure}{d}
\newcommand{\dx}{\dmesure\!}
\DeclareMathOperator{\Id}{Id}
\DeclareMathOperator{\Hom}{Hom}
\DeclareMathOperator{\cotanh}{cotanh}

\newcommand{\delr}{\partial_r}
\newcommand{\Jdelr}{J\!\delr}
\newcommand{\Ndr}{\nabla_{\delr}}
\newcommand{\Ldr}{\mathcal{L}_{\delr}}
\newcommand{\JE}{J\!E}

\renewcommand{\bar}{\overline}

\newcommand{\ALCH}{(\hyperlink{ALCH}{ALCH})~}
\newcommand{\AK}{(\hyperlink{AK}{AK})~}
\newcommand{\ALCHplus}{(\hyperlink{ALCHplus}{ALCH+})~}
\newcommand{\AKplus}{(\hyperlink{AKplus}{AK+})~}

\numberwithin{equation}{section}

\begin{document}


\begin{abstract}
In this article, we consider a complete, non-compact almost Hermitian manifold whose curvature is asymptotic to that of the complex hyperbolic plane.
Under natural geometric conditions, we show that such a manifold arises as the interior of a compact almost complex manifold whose boundary is a strictly pseudoconvex CR manifold.
Moreover, the geometric structure of the boundary can be recovered by analysing the expansion of the metric near infinity.
\end{abstract}

\maketitle




\section{Introduction}

\noindent The complex hyperbolic space is the unique simply connected, complete, Kähler manifold of constant negative holomorphic sectional curvature (we adopt the convention that this constant is $-1$).
It is the complex analogue of the real hyperbolic space, and similarly to its real counterpart, the complex hyperbolic space can be compactified by a sphere at infinity.
This sphere at infinity carries a natural geometric structure, which is closely related to the Riemannian geometry of the complex hyperbolic space: their respective groups of automorphisms are in one-to-one correspondence.
This structure is that of a strictly pseudoconvex CR manifold, namely, the CR sphere $(\mathbb{S},H,J)$.
If $\mathbb{S}$ is thought of as the unit sphere of $\C^N$, then $H = (T\mathbb{S})\cap (iT\mathbb{S})$ is the standard contact distribution, and $J$ is given by the multiplication by $i$ in $H$.
Set $\rho = e^{-r}$ with $r$ the distance function to a fixed point.
Then $\rho$ is a defining function for the boundary of the above compactification, and as $\rho \to 0$, the complex hyperbolic metric has the asymptotic expansion
\begin{equation}
\label{eq:complex_hyperbolic_metric}
\frac{1}{\rho^2}\dx \rho \otimes \dx \rho + \frac{1}{\rho^2}\theta\otimes\theta + \frac{1}{\rho}\gamma + o(1),
\end{equation}
with $\theta$ the standard contact form of $\mathbb{S}$, and $\gamma = \dx \theta|_{H\times H}(\cdot,J\cdot)$ the associated Levi-form.
The strict pseudoconvexity of the boundary means that the Levi-form is positive definite on $H$.

The aim of this paper is to construct a similar compactification by a strictly pseudoconvex CR structure for complete, non-compact, almost Hermitian manifolds satisfying some natural geometric conditions.
These conditions are the existence of a convex core (called an \emph{essential subset}), the convergence of the curvature tensor $R$ to that of the complex hyperbolic space $R^0$ near infinity, and the fact that the underlying almost complex structure $J$ is asymptotically Kähler at infinity.
More precisely, we show the following.

\begin{mainthm}
Let $(M,g,J)$ be a complete, non-compact, almost Hermitian manifold of real dimension at least $4$, which admits an essential subset.
Let $r$ be the distance function to any compact subset.
Assume that there exists $a > 1$ such that
\begin{equation}
\|R-R^0\|_g,\quad \|\nabla J\|_g,\quad \|\nabla R\|_g, \quad \text{and} \quad  \|\nabla^2 J\|_g = \mathcal{O}(e^{-ar}).
\end{equation}
Then $(M,J)$ is the interior of a compact almost complex manifold $(\bar{M},\bar{J})$, whose underlying almost complex structure $\bar{J}$ is continuous.
The hyperplane distribution $H_0 = (T\partial\bar{M})\cap (\bar{J}T\partial \bar{M})$ and the restriction $J_0 = \bar{J}|_{H_0}$ are of class $\mathcal{C}^1$.
Moreover, $H_0$ is a contact distribution, and $J_0$ is formally integrable, and $(\partial\bar{M},H_0,J_0)$ is a strictly pseudoconvex CR manifold.

\noindent In addition, the metric $g$ is asymptotically complex hyperbolic: there exists a defining function $\rho$ for the boundary, a $\mathcal{C}^1$ contact form $\eta^0$ calibrating $H_0$, and a continuous Carnot metric $\gamma$, with $\eta^0$ and $\gamma^0 = \gamma|_{H_0\times H_0} > 0$ of class $\mathcal{C}^1$, such that
\begin{equation}
\label{eq:CH_metric}
g \underset{\rho\to 0}{=} \frac{1}{\rho^2}\dx \rho \otimes \dx \rho + \frac{1}{\rho^2}\eta^0\otimes \eta^0 + \frac{1}{\rho} \gamma +
\begin{cases}
\mathcal{O}_g\left(\rho^{a-1}\right) & \text{if} \quad 1 < a < \frac{3}{2},\\
\mathcal{O}_g\left(\rho^{\frac{1}{2}}\ln \rho\right) & \text{if} \quad a = \frac{3}{2},\\
\mathcal{O}_g\left(\rho^{\frac{1}{2}}\right) & \text{if} \quad a > \frac{3}{2}.
\end{cases}
\end{equation}
The contact form and the Carnot metric are related by the relation $\dx\eta^0|_{H_0\times H_0}(\cdot,J_0\cdot) = \gamma^0$.
%
\end{mainthm}

This result gives a geometric characterisation of complete, non-compact, almost Hermitian manifolds admitting a compactification by a strictly pseudoconvex CR structure.
Notice the similarity between equations \eqref{eq:complex_hyperbolic_metric} and \eqref{eq:CH_metric}.
The real analogue of this result, involving a compactification by a conformal boundary for asymptotically locally real hyperbolic manifolds, has been proven by E. Bahuaud, J. M. Lee, T. Marsh and R. Gicquaud \cite{bahuaud_intrinsic_2008,bahuaud_conformal_2011,bahuaud_low_2018,bahuaud_holder_2008,gicquaud_conformal_2013}, pursuing the seminal work of M. T. Anderson and R. Schoen \cite{anderson_positive_1985}.

In a previous paper \cite{pinoy_asymptotic2021}, the author proved a similar result in the Kähler case.
The improvement here is twofold.
First, we are able to remove the Kähler assumption, which was of great importance in the previous proof.
Here, the almost complex structure is no more assumed to be parallel, and in fact, needs not even be formally integrable, nor the associated almost symplectic form needs to be closed.
In particular, the result applies to perturbations of asymptotically complex hyperbolic Kähler metrics which are only almost Hermitian.
Second, the strict pseudoconvexity of the boundary is obtained with an exponential decay of order $a > 1$, while the earlier version of this result needed a decay of order $a > \frac{3}{2}$.
Note that this has a cost: the Carnot metric can be shown to be $\mathcal{C}^1$ only in the direction of the contact distribution.
This is the reason why the extended almost complex structure $\bar{J}$ is only continuous in the transverse direction.
Both improvements imply that the set of examples to which the result applies is much increased.

A compactification by a CR structure for some complete, non-compact, Kähler manifolds was already given by J. Bland \cite{bland_existence_1985,bland_bounded_1989}, under assumptions that are rather analytic and not totally geometric.
To obtain a continuous compactification with no regularity on the CR structure, these assumptions imply the \emph{a posteriori} estimates $\|R-R^0\|_g, \|\nabla R\|_g = \mathcal{O}(e^{-4r})$\footnote{At first, one sees that these assumptions imply that $\|R-R^0\|_g = \mathcal{O}(e^{-3r})$ and $\|\nabla R\|_g = \mathcal{O}(e^{-4r})$.
Since on a Kähler manifold it holds that $\nabla R^0 = 0$, applying Kato's inequality to $R-R^0$ yields the claimed estimate.}.
A strictly pseudoconvex boundary of class $\mathcal{C}^1$ is obtained under assumptions that imply the even stronger estimates $\|R-R^0\|_g,\|\nabla R\|_g,\|\nabla^2 R\|_g = \mathcal{O}(e^{-5r})$.
It was proven by O. Biquard and M. Herzlich \cite{biquard_burns-epstein_2005} that for asymptotically complex hyperbolic Kähler-Einstein metrics in real dimension $4$, the curvature tensor has the form $R = R^0 + Ce^{-2r} + o_g(e^{-2r})$, where $C$ is a non-zero multiple of the Cartan tensor of the CR boundary.
It is known that the Cartan tensor vanishes exactly when the CR structure is locally equivalent to that of the sphere (such CR manifolds are called spherical).
Many examples are then not covered by J. Bland's results.

The paper is organized as follows.
In Section \ref{section:notations}, we set up the notations and explain the main idea of the proof of our main Theorem.
In Section \ref{section:metric_estimates}, we compute the expansion of the metric near infinity and prove the existence of the objects $\eta^0$ and $\gamma$, see Theorem \ref{theorem:metric_expansion}.
Section \ref{section:almost_complex} is dedicated to prove the existence of $J_0$, see Theorem \ref{theorem:almost_complex}.
At this step, $\eta^0$, $\gamma$ and $J_0$ are continuous tensor fields.
We show in Section \ref{section:higher_regularity} that they have higher regularity and that they induce a strictly pseudoconvex CR structure, see Theorems \ref{theorem:contact}, \ref{theorem:gamma_c^1} and \ref{theorem:strictly_pseudoconvex}.
Finally, we prove our main Theorem in Section \ref{section:compactification}.


\section{Preliminaries}
\label{section:notations}

\subsection{Notations}
Let $(M,g)$ be a Riemannian manifold.
Its Levi-Civita connection is denoted by $\nabla$.
Our convention on the Riemann curvature tensor is Besse's convention \cite{besse_einstein_2007}, namely
\begin{equation}
R(X,Y)Z = -\left(\nabla^2_{X,Y} Z - \nabla^2_{Y,X}Z\right) = \nabla_{[X,Y]}Z - \nabla_X(\nabla_YZ) + \nabla_Y(\nabla_XZ),
\end{equation}
for vector fields $X$, $Y$ and $Z$.
By abuse of notation, we still denote by $R$ its four times covariant version: this means that we write $R(X,Y,Z,T) = g(R(X,Y)Z,T)$ for vector fields $X$, $Y$, $Z$ and $T$.
With this convention, the sectional curvature of a tangent plane $P$ with orthonormal basis $\{u,v\}$ is $\sec(P) = \sec(u,v) = R(u,v,u,v)$.

\subsubsection*{Essential subsets and normal exponential map}
Following \cite{bahuaud_intrinsic_2008,bahuaud_conformal_2011,bahuaud_holder_2008,gicquaud_conformal_2013}, an \emph{essential subset} $K \subset M$ is a codimension $0$, compact, totally convex submanifold, with smooth boundary $\dK$ which is oriented by a unit outward vector field $\nu$, and such that $\sec(M\setminus K) < 0$.
In that case, the normal exponential map
\begin{equation}
\begin{array}{rccc}
\mathcal{E}\colon & \R_+ \times \dK & \longrightarrow & \bar{M\setminus K} \\
& (r,p) & \longmapsto & \exp_p(r\nu_p)
\end{array}
\end{equation}
is a diffeomorphism.
The level hypersurface at distance $r$ above $K$ is denoted by $\dK_r$.
For $r \geqslant 0$, $\mathcal{E}$ induces a diffeomorphism $\mathcal{E}_r\colon \dK \to \dK_r$ given by $\mathcal{E}_r(p)=\mathcal{E}(r,p)$; the induced Riemannian metric $\mathcal{E}_r^*g$ on $\dK$ is denoted by $g_r$.
Gauss Lemma states that $\mathcal{E}^*g = \dx r \otimes \dx r + g_r$.
Note that $g_0 = g|_{\dK}$.

The gradient of the distance function $r$ on $\bar{M\setminus K}$, called the \emph{radial vector field}, is denoted by $\delr$.
A \emph{radial geodesic} is a unit speed geodesic ray of the form $r \mapsto \mathcal{E}(r,p)$ with $p\in\dK$.
Note that the restriction of $\delr$ to a radial geodesic is its tangent vector field: therefore, $\delr$ satisfies the equation of geodesics $\Ndr \delr=0$.
More generally, a vector field $X$ on $\bar{M\setminus K}$ is called \emph{radially parallel} if $\Ndr X=0$.
The \emph{shape operator} $S$ is the field of symmetric endomorphisms on $\bar{M\setminus K}$ defined by $SX = \nabla_X\delr$.

The \emph{normal Jacobi field} on $\bar{M\setminus K}$ associated to a vector field $v$ on $\dK$ is defined by $Y_v = \mathcal{E}_*v$.
Such vector fields are orthogonal to and commute with the radial vector field $\delr$.
They satisfy the Jacobi field equation $\Ndr(\Ndr Y_v) = -R(\delr,Y_v)\delr$, and their restriction to any radial geodesic are thus Jacobi fields.
Normal Jacobi fields are related to the shape operator $S$ by the first order linear differential equation $\Ndr Y_v = SY_v$.

\subsubsection*{Almost Hermitian manifolds}
An \emph{almost Hermitian manifold} $(M,g,J)$ is a Riemannian manifold $(M,g)$ together with an almost complex structure $J$ which is compatible with the metric, in the sense that it induces linear isometries in the tangent spaces: one has $g(JX,JY) = g(X,Y)$ for all vector fields $X$ and $Y$.
Note that this implies that $J$ is skew-symmetric (in fact, these two properties are equivalent).
A tangent plane $P\subset TM$ is called \emph{$J$-holomorphic} (respectively \emph{totally real}) if $J\!P=P$ (respectively $J\!P\perp P$).
The \emph{constant $-1$ $J$-holomorphic sectional curvature tensor} $R^0$ on $(M,g,J)$ is defined by the equality
\begin{equation}
\label{eq:def_R0}
R^0(X,Y)Z = \frac{1}{4}\big( g(Y,Z)X - g(X,Z)Y + g(JY,Z)JX - g(JX,Z)JY  + 2g(X,JY)JZ\big)
\end{equation}
for $X$, $Y$ and $Z$ vector fields on $M$.
Similarly to the Riemann curvature tensor, we still denote by $R^0$ its fully covariant version, meaning that $R^0(X,Y,Z,T) = g(R^0(X,Y)Z,T)$ for all vector fields $X$, $Y$, $Z$ and $T$.
Note that $\|R^0\|_g \leqslant \frac{3}{2}$.
For any pair of orthogonal unit tangent vectors $u$ and $v$, $R^0(u,v,u,v) = -\frac{1}{4}(1+3g(Ju,v)^2)$;
the minimal value $-1$ (respectively the maximal value $-\frac{1}{4}$) is achieved precisely when $\{u,v\}$ spans a $J$-holomorphic plane (respectively a totally real plane).
In the specific case of the complex hyperbolic space, $R^0$ coincides with the curvature tensor of the complex hyperbolic metric (see \cite[Section IX.7]{kobayashi_foundations_1996-2}).

\subsubsection*{CR manifolds}
A \emph{CR manifold} (for Cauchy-Riemann) is a triplet $(M,H,J)$ where $H$ is a tangent distribution of hyperplanes and $J$ is an almost complex structure on $H$, such that the distribution $H^{1,0} = \{ X - iJX \mid X \in H\} \subset TM\otimes_{\R} \C$ is involutive (\emph{i.e.} $[X,Y]$ is a section of $H^{1,0}$ whenever $X$ and $Y$ are).
In this case, $J$ is said to be \emph{formally integrable}.
A CR manifold is called \emph{strictly pseudoconvex} if there exists a contact form $\eta$ calibrating the distribution $H$ (\emph{i.e.} $H=\ker \eta$ and $\dx \eta$ induces a non-degenerate $2$-form on $H$), and if the associated Levi form $\dx\eta|_{H\times H}(\cdot,J\cdot)$ is positive definite on $H$.

\subsection{The asymptotic conditions}
Throughout the paper, $(M,g,J)$ will denote a complete, non-compact, almost Hermitian manifold of dimension $2n+2\geqslant 4$, with an essential subset $K$.
We define the following asymptotic geometric conditions.

\begin{definition}[\ALCH and \AK conditions]
Let $(M,g,J)$ be a complete, non-compact, almost Hermitian manifold.
Let $r$ be the distance function to a compact subset.
\begin{enumerate}
\item We say that $(M,g,J)$ satisfies the \hypertarget{ALCH}{(ALCH)} condition of order $ a > 0$, for \emph{asymptotically locally complex hyperbolic}\footnote{For this condition implies that the local geometry at infinity resembles that of the complex hyperbolic space.}, if $\|R-R^0\|_g = \mathcal{O}(e^{-ar})$.

\item We say that $(M,g,J)$ satisfies the \hypertarget{AK}{(AK)} condition of order $a > 0$, for \emph{asymptotically Kähler}, if $\|\nabla J\|_g = \mathcal{O}(e^{-ar})$.
\end{enumerate}
\end{definition}

\begin{remark}
\label{remark:R_bounded}
Note that $\|R^0\|_g \leqslant \frac{3}{2}$.
The \ALCH condition of order $a > 0$ implies $\|R\|_g = \mathcal{O}(1)$.
\end{remark}

One readily verifies that the \ALCH condition implies that the sectional curvature of $M$ is bounded as follows: $-1 + \mathcal{O}(e^{-ar}) \leqslant \sec \leqslant - \frac{1}{4} + \mathcal{O}(e^{-ar})$.
The lower bound implies the following Lemma, proven in \cite[Proposition 2.5]{pinoy_asymptotic2021}.

\begin{lemma}
\label{lemma:S_bounded}
Assume that $(M,g,J)$ is a complete, non-compact, almost Hermitian manifold, admitting an essential subset $K$, and satisfying the \ALCH condition of order $a > 0$.
Let $S = \nabla \delr$ be the shape operator of the level hypersurfaces above $K$.
Then one has
\begin{equation}
\|S\|_g \leqslant 1 + \begin{cases}
\mathcal{O}\left(e^{-ar}\right) & \text{if} \quad 0 < a < 2,\\
\mathcal{O}\left((r+1)e^{-2r}\right) & \text{if} \quad a = 2,\\
\mathcal{O}\left(e^{-2r}\right) & \text{if} \quad a > 2.
\end{cases}
\end{equation}
In any case, one has $\|S\|_g = \mathcal{O}(1)$, and $\exp(\int_0^r \|S\|_g-1) = \mathcal{O}(1)$.
\end{lemma}

We also define the following analogous asymptotic conditions of higher order.

\begin{definition}[\ALCHplus and \AKplus conditions]
Let $(M,g,J)$ be a complete, non-compact, almost Hermitian manifold.
Let $r$ be the distance function to a compact subset.
\begin{enumerate}
\item We say that $(M,g,J)$ satisfies the \hypertarget{ALCHplus}{(ALCH+)} condition of order $ a > 0$ if one has the estimates $\|R-R^0\|_g = \mathcal{O}(e^{-ar})$ and $\|\nabla R\|_g = \mathcal{O}(e^{-ar})$.

\item We say that $(M,g,J)$ satisfies the \hypertarget{AKplus}{(AK+)} condition of order $a > 0$ if one has the estimates $\|\nabla J\|_g = \mathcal{O}(e^{-ar})$ and $\|\nabla^2 J\|_g = \mathcal{O}(e^{-ar})$.
\end{enumerate}
\end{definition}

\begin{remark}
\label{remark:nabla(R-R^0)}
Under the \AK condition of order $a > 0$, one has $\|\nabla R^0\|_g = \mathcal{O}(e^{-ar})$.
Thus, under the \AK condition of order $a > 0$, Kato's inequality shows that the \ALCHplus condition of order $a > 0$ is equivalent to $\|R-R^0\|_g \underset{r \to \infty}{\longrightarrow} 0$ and $\|\nabla(R-R^0)\|_g = \mathcal{O}(e^{-ar})$.
\end{remark}

In practice, $r$ will be the distance function to the essential subset $K$.
The constants involved in the previous estimates are global.
Moreover, in what follows, all estimates of the form $f = \mathcal{O}(h)$ will involve a constant that is global.
When built out of the choice of a reference frame (which will soon be called an
\emph{admissible frame}, see Definition \ref{definition:admissible}), the constant will be independent of that choice.
By the expressions $\|Y_u\|_g = \mathcal{O}(\|u\|_{g_0}e^r)$ or $Y_u = \mathcal{O}_g(\|u\|_{g_0}e^r)$, we mean that there exists $C > 0$ such that for any vector field $u$ on $\dK$, one has
\begin{equation}
\forall r \geqslant 0, \forall p \in \dK,\quad \|(Y_u)_{\mathcal{E}(r,p)}\|_g \leqslant C \|u_p\|_{g_0}e^r.
\end{equation}

\subsection{Outline of the proof}
If $(M,g,J)$ is assumed to be Kähler (that is, if $\nabla J=0$), the author showed in a previous paper \cite{pinoy_asymptotic2021} the following result.

\begin{thm*}[{\cite[Theorems A,B,C and D]{pinoy_asymptotic2021}}]
Let $(M,g,J)$ be a complete, non-compact, Kähler manifold admitting an essential subset $K$.
Assume that there is a constant $a>1$ such that the estimates $\|R-R^0\|_g,\|\nabla R\|_g=\mathcal{O}(e^{-ar})$ hold, where $r$ is the distance function to any compact subset.
Then on $\dK$, there exist a contact form $\eta$ of class $\mathcal{C}^1$, and a continuous symmetric positive bilinear form $\gamma$, positive definite on the contact distribution $H=\ker \eta$, such that
\begin{equation}
\label{eq:asymptotic_metric}
\mathcal{E}^*g = \dx r^2 + e^{2r}\eta\otimes \eta + e^r \gamma + \text{lower order terms}.
\end{equation}
If moreover $a>\frac{3}{2}$, then $\gamma$ is of class $\mathcal{C}^1$, and there exists a $\mathcal{C}^1$ formally integrable almost complex structure $J_H$ on $H$, such that $\gamma|_{H\times H} = \dx \eta(\cdot, J_H\cdot)$.
In particular, $(\dK,H,J_H)$ is a strictly pseudoconvex CR manifold.
\end{thm*}
\noindent Notice the similarity between equations \eqref{eq:CH_metric} and \eqref{eq:asymptotic_metric} by setting $\rho = e^{-r}$.
This result provides a compactification by a strictly pseudoconvex CR structure for a Kähler manifold whose curvature is asymptotically close to that of the complex hyperbolic space.
The proof is quite long, but can be summarised as follows:
\begin{enumerate}
\item For $\{J\nu,e_1,\ldots,e_{2n}\}$ an orthonormal frame on $\dK$, with $\nu$ the outward unit normal, let $\{\Jdelr,E_1,\ldots,E_{2n}\}$ denotes its parallel transport along radial geodesics.
For $r \geqslant 0$, define $\eta_r = \mathcal{E}_r^*(e^{-r}g(\cdot,\Jdelr))$, and $\eta^j_r = \mathcal{E}_r^*(e^{-\frac{r}{2}}g(\cdot,E_j))$, $j\in \{1,\ldots,2n\}$, which are local $1$-forms on $\dK$.

\item If $\|R-R^0\|_g = \mathcal{O}(e^{-ar})$, with $a > \frac{1}{2}$, then $\{\eta_r,\eta^1_r\ldots,\eta^{2n}_r\}_{r\geqslant 0}$ converges to continuous $1$-forms $\{\eta,\eta^1,\ldots,\eta^{2n}\}$.
This implies that the metric reads as in equation \eqref{eq:asymptotic_metric}, where $\gamma = \sum_{j=1}^{2n}\eta^j\otimes \eta^j$.
If moreover $a > 1$, volume comparison techniques show that the limit is a coframe.

\item If in addition, $\|\nabla R\|_g=\mathcal{O}(e^{-ar})$, then the family of $1$-forms $(\eta_r)_{r\geqslant 0}$ converges in $\mathcal{C}^1$ topology, the limit $\eta$ is of class $\mathcal{C}^1$, and is contact.
The proof uses several estimates, and tedious computations involving many curvature terms.

\item If $a>\frac{3}{2}$, then $(\eta_r^j)_{r\geqslant 0}$ locally uniformly converges in $\mathcal{C}^1$ topology, for any $j\in \{1,\ldots,2n\}$.
Hence, $\gamma$ is of class $\mathcal{C}^1$.

\item If $\varphi_r = \mathcal{E}_r^*\left(J - g(\cdot,\delr)\otimes \Jdelr) + g(\cdot,\Jdelr)\otimes \delr\right)$, then $(\varphi_r)_{r\geqslant 0}$ uniformly converges to a tensor $\varphi$ of class $\mathcal{C}^1$.
Its restriction to $H= \ker \eta$ gives the desired formally integrable almost complex structure $J_H$.
\end{enumerate}
The very first step of the proof crucially relies on the fact that $\Jdelr$ is parallel in the radial direction, and in fact, the equality $\nabla J = 0$ is used many times.
Note that the Kähler assumption is rather rigid: for instance, one has $\nabla J = 0$ if and only if the $2$-form $g(J\cdot,\cdot)$ is closed and $J$ is formally integrable.

In this paper, we extend and improve the results of \cite{pinoy_asymptotic2021}.
First, the Kähler condition is removed: in fact, neither the closedness of $g(J\cdot,\cdot)$ nor the formal integrability of $J$ need to be met.
We instead consider an almost Hermitian manifold $(M,g,J)$ whose almost complex structure $J$ is only parallel at infinity, by imposing the condition $\|\nabla^k J\|_g = \mathcal{O}(e^{-ar})$, $k\in\{1,2\}$.
Second, we show that the strict pseudoconvexity of the boundary can be obtained with $a > 1$ instead of $a > \frac{3}{2}$.
This sharper bound comes from deriving sharp geometric estimates in the direction of the contact structure.

In this context of this paper, the vector field $\Jdelr$ is not radially parallel, and one cannot even initiate the above strategy as it stands.
The main trick is to prove the existence, under our assumptions, of a unit vector field $E_0$ on $\bar{M \setminus K}$ that is radially parallel, and that satisfies $\|E_0-\Jdelr\|_g = \mathcal{O}(e^{-ar})$.
This latter vector field is unique.
One can then consider a reference frame $\{E_0,\ldots,E_{2n}\}$ having nice properties, which we call an \emph{admissible frame} (see Definition \ref{definition:admissible} below), and try to mimic the above proof.
The counterpart is that the computations become longer and more involved; one also needs to show numerous extra estimates.


\section{Metric estimates}
\label{section:metric_estimates}

\noindent This section is dedicated to the derivation of the expansion near infinity of the metric $g$ under the \ALCH and \AK conditions.
We first define the notion of admissible frames, which simplify future computations.
We then derive estimates on the asymptotic expansion of normal Jacobi fields, which turns out to be the main ingredients to show our results.

\subsection{Admissible frames}
\label{subsection:admissible_frame}
We give a construction for some parallel orthonormal frames along radial geodesics in which later computations will be easier.
For $v$ a vector field on $\dK$, let $V$ be the vector field on $\bar{M\setminus K}$ obtained by the parallel transport of $v$ along radial geodesics.
Finally, for $r \geqslant 0$, define $\beta_r(v) = g(\Jdelr,V)|_{\dK_r}$.
This defines a family of $1$-forms $(\beta_r)_{r\geqslant 0}$ on $\dK$.

\begin{lemma}
\label{lemma:beta_r_converges}
Let $(M,g,J)$ be a complete, non-compact, almost Hermitian manifold of dimension at least $4$, with essential subset $K$.
Assume that it satisfies the \AK condition of order $a > 0$.
Then there exists a continuous $1$-form  $\beta$ on $\dK$ such that
\begin{equation}
\label{eq:beta_r-beta}
\beta_r - \beta = \mathcal{O}_{g_0}(e^{-ar}).
\end{equation}
\end{lemma}

\begin{proof}
Fix $v$ a vector field on $\dK$ and $r \geqslant 0$.
Both $\delr$ and $V$ are radially parallel, so that one has $\beta_r(v)-\beta_0(v) = \int_0^r \delr  g(\Jdelr,V) = \int_0^r g((\Ndr J)\delr,V)$.
By the \AK assumption, there exists $C > 0$ such that $\|\nabla J\|_g \leqslant Ce^{-ar}$.
The Cauchy-Schwarz inequality now implies that $\int_0^r\|g((\Ndr J)\delr, V)\| \leqslant \int_0^r \|\nabla J\|_g \|V\|_{g} \leqslant C\frac{1-e^{-ar}}{a} \|v\|_{g_0}$.
Therefore, $(\beta_r(v))_{r\geqslant 0}$ pointwise converges: let $\beta(v)$ to be its pointwise limit.
It defines a pointwise linear form on the tangent spaces of $\dK$, satisfying
\begin{equation}
\label{eq:beta-beta_r}
\left|\beta(v)-\beta_r(v)\right|
= \left| \int_r^{\infty} g((\Ndr J)\delr,V) \right|
\leqslant \int_r^{\infty} \left|g((\Ndr J)\delr,V)\right|
\leqslant \frac{C}{a}e^{-ar}\|v\|_{g_0},
\end{equation}
from which is derived equation \eqref{eq:beta_r-beta}.
The convergence is thus uniform, and $\beta$ is continuous.

We shall now show that $\beta$ is nowhere vanishing.
For all $r \geqslant 0$, one has $\|\beta_r\|_{g_0} = 1$ pointwise.
Indeed, for any $v$, Cauchy-Schwarz inequality implies that $|\beta_r(v)| \leqslant \|V\|_g  = \|v\|_{g_0}$.
Equality is reached for $v = \iota_r^{-1}(\Jdelr)$, where $\iota_r\colon T\dK \to T\dK_r$ is induced by the parallel transport along radial geodesics.
It follows that $\|\beta\|_{g_0} = 1$ pointwise, and that $\beta$ is nowhere vanishing.
\end{proof}

\begin{definition}
\label{definition:admissible}
Let $(M,g,J)$ be a complete, non-compact, almost Hermitian manifold of dimension at least $4$, with essential subset $K$.
Assume that it satisfies the \AK condition of order $a > 0$.
Let $U\subset \dK$ be an open subset on which the continuous distribution $\ker \beta$ is trivialisable.
Let $\{e_0,\ldots,e_{2n}\}$ be an orthonormal frame on $U$ such that $\beta(e_0) > 0$ and $\beta(e_j) = 0$ if $j\in \{1,\ldots,2n\}$.
The associated \emph{admissible frame} $\{E_0,\ldots,E_{2n}\}$ on the cone $E(\R_+\times U)$ is defined as the parallel transport of $\{e_0,\ldots,e_{2n}\}$ along the radial geodesics.
\end{definition}

If $\{E_0,\ldots,E_{2n}\}$ is an admissible frame, then $\{\delr,E_0,\ldots,E_{2n}\}$ is an orthonormal frame on the cone $E(\R_+\times U)$ whose elements are parallel in the radial direction even though they need not be differentiable in the directions that are orthogonal to $\delr$.
In the following, we will often refer to admissible frames without mentioning the open subset $U\subset \dK$ on which they are defined.

\begin{lemma}
\label{lemma:beta(e_0)=1}
Let $(M,g,J)$ be a complete, non-compact, almost Hermitian manifold of dimension at least $4$, with essential subset $K$.
Assume that it satisfies the \AK condition of order $a > 0$.
Let $\{E_0,\ldots,E_{2n}\}$ be an admissible frame.
Then $\beta(e_0) = 1$.
\end{lemma}

\begin{proof}
One has $1 = \|\Jdelr \|_g^2 = \sum_{j=0}^{2n} \beta_r(e_j)^2$.
The result follows by taking the limit as $r \to \infty$.
\end{proof}

\begin{corollary}
\label{corollary:frame_estimates}
Let $(M,g,J)$ be a complete, non-compact, almost Hermitian manifold of dimension at least $4$, with essential subset $K$.
Assume that it satisfies the \AK condition of order $a > 0$.
Let $\{E_0,\ldots,E_{2n}\}$ be an admissible frame and $\delta$ be the Kronecker symbol.
Then
\begin{enumerate}
\item $g(\Jdelr,E_j) - \delta_{0j} = \mathcal{O}(e^{-ar})$ for $j\in \{0,\ldots,2n\}$,
\item $E_0 - \Jdelr = \mathcal{O}_g(e^{-ar})$.
\end{enumerate}
\end{corollary}

\begin{proof}
The first point is a consequence of the equality $g(\Jdelr,E_j)=\beta_r(e_j)$ and of equation \eqref{eq:beta-beta_r}.
For the second point, notice that
\begin{equation}
\label{eq:E_0-Jdelr}
E_0-\Jdelr = \sum_{j=0}^{2n}g(E_0-\Jdelr,E_j)E_j = \sum_{j=0}^{2n}(\delta_{0j}- g(\Jdelr,E_j))E_j,
\end{equation}
from which is derived the claimed estimate.
\end{proof}

\begin{remark}
\label{remark:kahler1}
One easily shows that the vector field $E_0$ is the unique unit vector field $X$ on $E(\R_+\times U)$ such that $\Ndr X = 0$ and $g(X,\Jdelr) = 1 + o(1)$.
If $(M,g,J)$ is Kähler (if $\nabla J = 0$), then $\Ndr\Jdelr = 0$, and thus $E_0 = \Jdelr$.
In this specific case, admissible frames can be chosen to be smooth, and correspond to the \emph{radially parallel orthonormal frames} defined in \cite{pinoy_asymptotic2021}.
\end{remark}

\begin{proposition}
\label{proposition:curvature_estimates}
Let $(M,g,J)$ be a complete, non-compact, almost Hermitian manifold of dimension at least $4$, with essential subset $K$.
Assume that it satisfies the \ALCH and \AK conditions of order $a > 0$.
Let $\{E_0,\ldots,E_{2n}\}$ be an admissible frame.
Then
\begin{enumerate}
\item $\sec(\delr,E_0) + 1 = \mathcal{O}(e^{-ar})$,
\item $\sec(\delr,E_j) + \frac{1}{4} = \mathcal{O}(e^{-ar})$ for $j \in \{1,\ldots,2n\}$,
\item $R(\delr,E_i,\delr,E_k) = \mathcal{O}(e^{-ar})$ for any $i \neq j \in \{0,\ldots,2n\}$.
\end{enumerate}
\end{proposition}

\begin{proof}
We prove the first point, the other being shown similarly.
One readily verifies from the definition of $R^0$ that $R^0(\delr,\Jdelr,\delr,\Jdelr) = -1$, and therefore, it holds that
\begin{equation}
\begin{split}
\sec(\delr,E_0)
&= R^0(\delr,\Jdelr + (E_0-\Jdelr), \delr, \Jdelr + (E_0-\Jdelr))+ (R-R^0)(\delr,E_0,\delr,E_0) \\
&= -1 + 2R^0(\delr,E_0-\Jdelr,E_0,\Jdelr)
+ R^0(\delr,E_0-\Jdelr,\delr,E_0-\Jdelr)\\
&\quad + (R-R^0)(\delr,E_0,\delr,E_0).
\end{split}
\end{equation}
The definition of $R^0$ (see equation \eqref{eq:def_R0}) yields $\|R^0\|_g \leqslant \frac{3}{2}$, and the result follows from the \ALCH assumption and from the second point of Corollary \ref{corollary:frame_estimates}.
\end{proof}

\subsection{Associated coframes and normal Jacobi fields estimates}
\label{subsection:associated_coframes}
Recall that for $r \geqslant 0$, the diffeomorphism $\mathcal{E}_r\colon \dK \to \dK_r$ is defined by $\mathcal{E}_r(p) = \mathcal{E}(r,p)$.
\begin{definition}
Let $(M,g,J)$ be a complete, non-compact, almost Hermitian manifold with essential subset $K$.
Assume that it satisfies the \AK condition of order $a > 0$.
Let $\{E_0,\ldots,E_{2n}\}$ be an admissible frame on the cone $E(\R_+\times U)$.
The \emph{associated coframe} $\{\eta^0_r,\ldots,\eta^{2n}_r\}_{r \geqslant 0}$ on $U$ is defined by
\begin{equation}
\eta^0_r =  \mathcal{E}_r^* \left(e^{-r} g(\cdot,E_0)\right) \quad \text{and}
\quad \eta^j_r = \mathcal{E}_r^*\left(e^{-\frac{r}{2}}g(\cdot,E_j)\right)
\quad \text{if} \quad j\in \{1,\ldots,2n\}.
\end{equation}
\end{definition}

In any admissible frame, the normal Jacobi field $Y_v$ associated to the vector field $v$ on $\dK$ reads
\begin{equation}
\label{eq:expression_jacobi_field}
Y_v = \eta^0_r(v) e^r E_0
+ \sum_{j=1}^{2n} \eta^j_r(v) e^{\frac{r}{2}}E_j.
\end{equation}
Applying twice the differential operator $\Ndr$ to this last equality, one has
\begin{equation}
\begin{split}
\Ndr(\Ndr Y_v) &=  \left(\delr^2 \eta^0_r(v)+ 2\delr \eta^0_r(v) + \eta^0_r(v) \right)e^r E_0 \\
&\quad + \sum_{j=1}^{2n}\left(\delr^2\eta^j_r(v) + \delr\eta^j_r(v) + \frac{1}{4}\eta^j_r(v) \right)e^{\frac{r}{2}}E_j.
\end{split}
\end{equation}
Recall that radial Jacobi fields are actual Jacobi fields, which means that they satisfy the second order linear differential equation $\Ndr(\Ndr Y_v) = -R(\delr,Y_v)\delr$.
An identification of the components of $\Ndr(\Ndr Y_v)$ in the given admissible frame shows that the coefficients $\{\eta^j_r(v)\}_{j \in \{0,\ldots,2n\}}$ satisfy the differential system
\begin{equation}
\begin{cases}
\displaystyle \delr^2\eta^0_r(v) + 2 \delr \eta^0_r(v) =
\sum_{k=0}^{2n} u^0_k \eta^k_r(v), \\
\displaystyle \delr^2\eta^j_r(v) + \phantom{2} \delr \eta^j_r(v) =
\sum_{k=0}^{2n} u^j_k \eta^k_r(v), & j\in \{1,\ldots,2n\},
\end{cases}
\end{equation}
where the functions $\{u^j_k\}_{j,k\in \{0,\ldots,2n\}}$ are defined by
\begin{equation}
u^j_k = -
\begin{cases}
\sec(\delr,E_0) + 1 & \text{if} \quad j=k=0,\\
e^{-\frac{r}{2}}R(\delr,E_0,\delr,E_k) & \text{if} \quad j=0, k\neq 0, \\
e^{\frac{r}{2}} R(\delr,E_k,\delr,E_0) & \text{if} \quad j\neq 0, k=0, \\
R(\delr,E_j,\delr,E_k) & \text{if} \quad j,k \in \{1,\ldots,2n\}, j\neq k,\\
\sec(\delr,E_j) + \frac{1}{4} & \text{if} \quad j,k\in\{1,\ldots,2n\}, j=k.
\end{cases}
\end{equation}
Proposition \ref{proposition:curvature_estimates} implies that one has the uniform estimates $|u^j_k| = \mathcal{O}(e^{-(a-\frac{1}{2})r})$.
Combining the proofs of \cite[Propositions 3.7 \& 3.14]{pinoy_asymptotic2021}, relying on successive integrations, an application of Grönwall's Lemma, and a bootstrap argument, one obtains the following result.
The last claim relies on estimates on the growth of the volume (see \cite[Propositions 2.7 \& 3.13]{pinoy_asymptotic2021}).

\begin{proposition}
\label{proposition:eta_r_convergence}
Let $(M,g,J)$ be a complete, non-compact, almost Hermitian manifold of dimension at least $4$, with essential subset $K$.
Assume that it satisfies the \ALCH and \AK conditions of order $a>\frac{1}{2}$.
Let $\{\eta^0_r,\ldots,\eta^{2n}_r\}_{r \geqslant 0}$ be the coframes associated to an admissible frame on $U\subset \dK$.
Then there exists continuous $1$-forms $\{\eta^0,\ldots,\eta^{2n}\}$ on $U$
\begin{equation}
\begin{split}
\partial_r \eta^0_r, \quad \eta^0_r - \eta^0 &=
\begin{cases}
\mathcal{O}_{g_0}\left(e^{-ar}\right) & \text{if} \quad
\frac{1}{2} < a <\frac{3}{2}, \\
\mathcal{O}_{g_0}\left((r+1)e^{-\frac{3}{2}r}\right) & \text{if} \quad
a = \frac{3}{2}, \\
\mathcal{O}_{g_0}\left(e^{-\frac{3}{2}r}\right) & \text{if} \quad
a > \frac{3}{2},
\end{cases}\\
\forall j \in \{1,\ldots,2n\},\quad
\partial_r \eta^j_r, \quad \eta^j_r - \eta^j &=
\begin{cases}
\mathcal{O}_{g_0}\left(e^{-(a-\frac{1}{2})r}\right) & \text{if} \quad
\frac{1}{2} < a <\frac{3}{2}, \\
\mathcal{O}_{g_0}\left((r+1)e^{-r}\right) & \text{if} \quad
a = \frac{3}{2}, \\
\mathcal{O}_{g_0}\left(e^{-r}\right) & \text{if} \quad
a > \frac{3}{2}.
\end{cases}
\end{split}
\end{equation}
If furthermore one assumes that $a > 1$, the family $\{\eta^0,\ldots,\eta^{2n}\}$ is a continuous coframe on $U$.
\end{proposition}

\begin{corollary}
\label{corollary:eta_r_uniform_bound}
If $a > \frac{1}{2}$, then $\|\eta^j_r\|_{g_0}$ is bounded independently of $r$, $j$, the choice of an admissible frame, and $U$.
\end{corollary}

\begin{proof}
For $j\in \{0,\ldots, 2n\}$ and $r \geqslant 0$, write $\eta^j_r = \eta^j_0 + \int_0^r \delr \eta^j_r$.
Notice that $\|\eta^j_0\|_{g_0} = 1$.
Then by Proposition \ref{proposition:eta_r_convergence}, $\|\eta^j_r\|_{g_0}  \leqslant \|\eta^j_0\|_{g_0} + \int_0^r \|\delr \eta^j_r\|_{g_0} \leqslant 1 + \int_0^{\infty} \|\delr \eta^j_r\|_{g_0} = \mathcal{O}(1)$.
\end{proof}

Recall that a normal Jacobi field $Y_v$ satisfies $\Ndr Y_v = SY_v$.
The following corollary is an immediate consequence of Proposition \ref{proposition:eta_r_convergence}.

\begin{corollary}
\label{corollary:Y_v_asymptotic}
In any admissible frame, the normal Jacobi field $Y_v$ associated to a vector field $v$ on $\dK$ satisfies
\begin{equation}
\label{eq:Y_v_asymptotic}
Y_v = \eta^0(v) e^r E_0 + \sum_{j=1}^{2n} \eta^j(v)e^{\frac{r}{2}} E_j +
\begin{cases}
\mathcal{O}_g\left(\|v\|_{g_0} e^{-(a-1)r}\right) & \text{if} \quad
\frac{1}{2} < a <\frac{3}{2},\\
\mathcal{O}_g\left(\|v\|_{g_0} (r+1)e^{-\frac{r}{2}}\right) & \text{if} \quad
a = \frac{3}{2},\\
\mathcal{O}_g\left(\|v\|_{g_0} e^{-\frac{r}{2}}\right) & \text{if} \quad
a > \frac{3}{2},\\
\end{cases}
\end{equation}
and
\begin{equation}
\label{eq:SY_v_asymptotic}
SY_v = \eta^0(v) e^r E_0 + \sum_{j=1}^{2n} \frac{1}{2}\eta^j(v)e^{\frac{r}{2}} E_j +
\begin{cases}
\mathcal{O}_g\left(\|v\|_{g_0} e^{-(a-1)r}\right) & \text{if} \quad
\frac{1}{2} < a <\frac{3}{2},\\
\mathcal{O}_g\left(\|v\|_{g_0} (r+1)e^{-\frac{r}{2}}\right) & \text{if} \quad
a = \frac{3}{2},\\
\mathcal{O}_g\left(\|v\|_{g_0} e^{-\frac{r}{2}}\right) & \text{if} \quad
a > \frac{3}{2}.\\
\end{cases}
\end{equation}
As a consequence, one has the global estimates $Y_v, SY_v = \mathcal{O}_g(\|v\|_{g_0}e^r)$.
If moreover, $v$ is everywhere tangent to $\ker \eta^0$, then $Y_v, SY_v = \mathcal{O}_g(\|v\|_{g_0}e^{\frac{r}{2}})$.
\end{corollary}

\begin{remark}
Note that although the estimates of Proposition \ref{proposition:eta_r_convergence} are not uniform in all directions, they contribute equally to the lower order term in equations \eqref{eq:Y_v_asymptotic} and \eqref{eq:SY_v_asymptotic} thanks to the remaining exponential factors.
\end{remark}

\subsection{Global consequences and metric estimates}
\label{subsection:metric estimates}
We shall now highlight global consequences of the study conducted in Subsections \ref{subsection:admissible_frame} and \ref{subsection:associated_coframes}.
We then prove the first of our main results.

\begin{lemma}
\label{lemma:e_0_global}
Assume that $(M,g,J)$ satisfies the \AK condition of order $a > 0$.
Then the local vector field $e_0$ defined in Definition \ref{definition:admissible} defines a global continuous vector field on $\dK$, independently of the construction of any admissible frame.
\end{lemma}

\begin{proof}
The $1$-form $\beta$ defined in Lemma \ref{lemma:beta_r_converges} is continuous and nowhere vanishing.
Hence, the distribution $\ker \beta \subset T\dK$ is a continuous distribution of hyperplanes.
It follows that its $g_0$-orthogonal complement $L$ is a well-defined and continuous line bundle.
Notice that the restriction of $\beta$ trivialises $L$.
It follows that $e_0$ is the unique section of $L$ that is positive for $\beta$, and of unit $g_0$-norm.
This concludes the proof.
\end{proof}

The family of $1$-forms $\{\eta^0_r\}_{r \geqslant 0}$ is then globally defined on $\dK$, independently of the choice of the admissible frame.
As a consequence, one has the following global version of Proposition \ref{proposition:eta_r_convergence}.

\begin{proposition}
\label{proposition:eta^0_global}
Let $(M,g,J)$ be a complete, non-compact, almost Hermitian manifold of dimension at least $4$, admitting an essential subset $K$.
Assume that it satisfies the \ALCH and \AK condition of order $a > \frac{1}{2}$.
Then there exists a continuous $1$-form $\eta^0$ on $\dK$ such that
\begin{equation}
\begin{split}
\partial_r \eta^0_r,\quad  \eta^0_r - \eta^0 &=
\begin{cases}
\mathcal{O}_{g_0}\left(e^{-ar}\right) & \text{if} \quad
\frac{1}{2} < a <\frac{3}{2}, \\
\mathcal{O}_{g_0}\left((r+1)e^{-\frac{3}{2}r}\right) & \text{if} \quad
a = \frac{3}{2}, \\
\mathcal{O}_{g_0}\left(e^{-\frac{3}{2}r}\right) & \text{if} \quad
a > \frac{3}{2}.
\end{cases}
\end{split}
\end{equation}
If furthermore one assumes that $a > 1$, then $\eta^0$ is nowhere vanishing.
\end{proposition}

The following Corollary is a straightforward application of the triangle inequality and of Corollary \ref{corollary:eta_r_uniform_bound}.

\begin{corollary}
\label{corollary:eta_tensor_eta}
One has the following estimates
\begin{equation}
\eta^0_r \otimes \eta^0_r - \eta^0 \otimes \eta^0 =
\begin{cases}
\mathcal{O}_{g_0}\left(e^{-ar} \right) & \text{if} \quad \frac{1}{2} < a < \frac{3}{2},\\
\mathcal{O}_{g_0}\left((r+1)e^{-\frac{3}{2}r} \right) & \text{if} \quad a = \frac{3}{2},\\
\mathcal{O}_{g_0}\left(e^{-\frac{3}{2}r} \right) & \text{if} \quad a > \frac{3}{2}.
\end{cases}
\end{equation}
\end{corollary}


From Gauss's Lemma, the Riemannian metric $g$ reads as $\mathcal{E}^*g = \dx r \otimes \dx r + g_r$, with $(g_r)_{r \geqslant 0}$ the smooth family of Riemannian metrics on $\dK$ defined by $g_r = \mathcal{E}_r^* g$.
By construction, the first term that appears in the asymptotic expansion of the metric $g$ near infinity is $e^{2r}\eta^0 \otimes \eta^0$.

\begin{definition}
For $r\geqslant 0$, $\gamma_r$ is defined as $\gamma_r = e^{-r}( g_r - e^{2r} \eta^0_r \otimes \eta^0_r)$.
\end{definition}

By definition, $(\gamma_r)_{r\geqslant 0}$ is a family of symmetric $2$-tensors on $\dK$.
Let $\{\eta^0_r,\ldots,\eta^{2n}_r\}_{r \geqslant 0}$ be the coframes associated to an admissible frame $\{E_0,\ldots,E_{2n}\}$.
Then locally, $\gamma _r = \sum_{j=1}^{2n} \eta^j_r\otimes \eta^j_r$.
Consequently, $\gamma_r$ is positive semi-definite, and is positive definite on $\ker \eta^0_r$, for any $r \geqslant 0$.
The following proposition shows that $(\gamma_r)_{r \geqslant 0}$ converges to some tensor that shares similar properties.

\begin{proposition}
\label{proposition:gamma_r_convergence}
Let $(M,g,J)$ be a complete, non-compact, almost Hermitian manifold of dimension at least $4$, and admitting an essential subset $K$. Assume that it satisfies the \ALCH and \AK conditions of order $a > \frac{1}{2}$.
Then there exists a continuous positive semi-definite symmetric $2$-tensor $\gamma$ on $\dK$, which we call the Carnot metric, such that
\begin{equation}
\label{eq:gamma_r-gamma}
\gamma_r - \gamma =
\begin{cases}
\mathcal{O}_{g_0}\left(e^{-(a-\frac{1}{2})r}\right) & \text{if} \quad \frac{1}{2} < a < \frac{3}{2},\\
\mathcal{O}_{g_0}\left((r+1)e^{-r}\right) & \text{if} \quad a = \frac{3}{2},\\
\mathcal{O}_{g_0}\left(e^{-r}\right) & \text{if} \quad a > \frac{3}{2}.
\end{cases}
\end{equation}
If furthermore one assumes that $a > 1$, then $\gamma$ is positive definite on $\ker \eta^0$.
\end{proposition}

\begin{proof}
For $r \geqslant 0$, one has $g_r = e^{2r}\eta^0_r\otimes \eta^0_r + e^r \gamma_r$.
Let $\{\eta^0_r,\ldots,\eta^{2n}\}_{r \geqslant 0}$ be the coframes associated with an admissible frame.
Locally, one can express $\gamma_r$ as $\gamma_r = \sum_{j=1}^{2n} \eta^j_r\otimes \eta^j_r$.
Therefore, $(\gamma_r)_{r \geqslant 0}$ converges pointwise to a limit we call $\gamma$ which is locally given by $\sum_{j=1}^{2n} \eta^j\otimes \eta^j$.
In addition, one has the local expression
$\gamma_r - \gamma = \sum_{j=1}^{2n} \eta^j_r\otimes (\eta^j_r-\eta^j) + (\eta^j_r-\eta^j) \otimes \eta^j$.
The global estimates \eqref{eq:gamma_r-gamma} now follow from the triangle inequality and from an application of Proposition \ref{proposition:eta_r_convergence} and Corollary \ref{corollary:eta_r_uniform_bound}.
As a consequence, $\gamma$ is a continuous symmetric positive semi-definite $2$-tensor.
If $a > 1$, then $\{\eta^0,\ldots,\eta^{2n}\}$ is a coframe (Proposition \ref{proposition:eta_r_convergence}), and $\gamma$ is hence positive definite on $\ker \eta^0$.
\end{proof}

The previous study implies the following comparison between quadratic forms.

\begin{corollary}
\label{corollary:pinching_g_g_r}
If $a > 1$, there exists a constant $\lambda > 1$ such that for all $r \geqslant 0$, the comparison between quadratic forms $\frac{1}{\lambda} e^{r}g_0 \leqslant g_r \leqslant \lambda e^{2r} g_0$ holds.
\end{corollary}

\begin{proof}
For $r \geqslant 0$, $\eta^0_r \otimes \eta^0_r$ and $\gamma_r$ are positive symmetric $2$-tensors.
Define $q_r = \eta_r^0\otimes \eta_r^0 + \gamma_r$, which is a Riemannian metric on $\dK$.
From $g_r = e^{2r}\eta^0_r \otimes \eta^0_r + e^r \gamma_r$, one readily checks that
\begin{equation}
\label{eq:pinching_q_r}
\forall r \geqslant 0,\quad e^r q_r \leqslant g_r \leqslant e^{2r}q_r.
\end{equation}
According to Propositions \ref{proposition:eta^0_global} and \ref{proposition:gamma_r_convergence}, $q_r$ uniformly converges to the continuous Riemannian metric $q_{\infty} = \eta^0 \otimes \eta^0 + \gamma$ as $r\to \infty$.
Let $S^{g_0}\dK$ be the unit sphere bundle of $(\dK,g_0)$, which is compact by compactness of $\dK$.
The map $(r,v) \in [0,\infty]\times S^{g_0}\dK \mapsto q_r(v,v)\in (0,\infty)$ is then continuous on the compact space $[0,\infty]\times S^{g_0}\dK$.
Therefore, there exists $\lambda > 1$ such that for all $r\geqslant 0$,  $\frac{1}{\lambda} \leqslant q_r \leqslant \lambda$ on $S^{g_0}\dK$.
The result now follows from equation \eqref{eq:pinching_q_r} and from the homogeneity of quadratic forms.
\end{proof}

We shall now show the first of our main results.

\begin{thm}
\label{theorem:metric_expansion}
Let $(M,g,J)$ be a complete, non-compact, almost Hermitian manifold of dimension at least $4$, with essential subset $K$.
Assume that it satisfies the \ALCH and \AK assumptions of order $a > \frac{1}{2}$.
Then on $\dK$, there exists a continuous $1$-form $\eta^0$ and a continuous positive semi-definite symmetric $2$-tensor $\gamma$, such that in the normal exponential map $E$, the Riemannian metric $g$ reads
\begin{equation}
\label{eq:metric_estimates}
\mathcal{E}^*g = \dx r \otimes \dx r + e^{2r} \eta^0 \otimes \eta^0
+ e^r \gamma +
\begin{cases}
\mathcal{O}_{g_0}\left(e^{(2-a)r}\right)
& \text{if} \quad \frac{1}{2} < a < \frac{3}{2}, \\
\mathcal{O}_{g_0}\left((r+1)e^{\frac{r}{2}}\right)
& \text{if} \quad a = \frac{3}{2}, \\
\mathcal{O}_{g_0}\left(e^{\frac{r}{2}}\right)
& \text{if} \quad a > \frac{3}{2}.
\end{cases}
\end{equation}
If furthermore one assumes that $a > 1$, then $\eta^0$ is nowhere vanishing, and $\gamma$ is positive definite on the distribution of hyperplanes $\ker \eta^0$.
\end{thm}

\begin{proof}
Let $(\eta^0_r)_{r \geqslant 0}$, $(\gamma_r)_{r \geqslant 0}$ and their limits $\eta^0$ and $\gamma$ be given by
Propositions \ref{proposition:eta^0_global} and \ref{proposition:gamma_r_convergence}.
By construction, one has
\begin{equation}
\mathcal{E}^*g = \dx r \otimes \dx r + e^{2r}\eta^0_r \otimes \eta^0_r + e^r \gamma_r
= \dx r \otimes \dx r + e^{2r}\eta^0 \otimes \eta^0 + e^r \gamma + \varepsilon_r,
\end{equation}
with $\varepsilon_r = e^{2r}\left(\eta^0_r \otimes \eta^0_r - \eta^0 \otimes \eta^0\right) + e^r (\gamma_r - \gamma)$.
Estimates \eqref{eq:metric_estimates} now follow from Corollary \ref{corollary:eta_tensor_eta} (estimates on $\eta^0_r\otimes\eta^0_r - \eta^0\otimes\eta^0$)
and Proposition \ref{proposition:gamma_r_convergence} (estimates on $\gamma_r-\gamma$).
Ultimately, if $a > 1$, the last claim follows from Propositions \ref{proposition:eta^0_global} ($\eta^0$ is nowhere vanishing) and \ref{proposition:gamma_r_convergence} ($\gamma$ is positive semi-definite, positive definite on $\ker \eta^0$).
\end{proof}

\begin{remark}
\label{remark:g-ghat}
Setting $\widehat{g} = \mathcal{E}_*(\dx r\otimes \dx r + e^{2r} \eta^0\otimes \eta^0 + e^r \gamma)$ on $\bar{M\setminus K}$, Corollary \ref{corollary:pinching_g_g_r} shows that estimates \eqref{eq:metric_estimates} read
\begin{equation}
g - \widehat{g} =
\begin{cases}
\mathcal{O}_{g}\left(e^{-(a-1)r}\right)
& \text{if} \quad \frac{1}{2} < a < \frac{3}{2}, \\
\mathcal{O}_{g}\left((r+1)e^{-\frac{r}{2}}\right)
& \text{if} \quad a = \frac{3}{2}, \\
\mathcal{O}_{g}\left(e^{-\frac{r}{2}}\right)
& \text{if} \quad a > \frac{3}{2}.
\end{cases}
\end{equation}
If $\eta^0$ were a contact form and $\gamma$ a Carnot metric on its kernel distribution, then $g$ would be asymptotically complex hyperbolic in the sense of \cite{biquard_metriques_2000,biquard_burns-epstein_2005}.
\end{remark}

\subsection{Estimates on the shape operator}
Before we conclude this section, we give another consequence of the previous study: we derive asymptotic estimates on the shape operator $S$.
First, we introduce a natural vector field $\xi_0$, which is closely related to $S$.

\begin{definition}
The vector fields $(\xi_0^r)_{r \geqslant 0}$ on $\dK$ are defined as $\xi_0^r = \mathcal{E}_r^* (e^r E_0)$.
\end{definition}

\begin{proposition}
\label{proposition:xi_0_estimates}
Let $(M,g,J)$ be a complete, non-compact, almost Hermitian manifold of dimension at least $4$, admitting an essential subset $K$.
Assume that it satisfies the \ALCH and \AK conditions of order $a > 1$.
Then there exists a continuous vector field $\xi_0$ on $\dK$ such that
\begin{equation}
\label{eq:estimates_xi_0}
\xi_0^r - \xi_0\ =
\begin{cases}
\mathcal{O}_{g_0}\left(e^{-(a-\frac{1}{2})r}\right) & \text{if} \quad
1 < a <\frac{3}{2}, \\
\mathcal{O}_{g_0}\left((r+1)e^{-r}\right) & \text{if} \quad
a = \frac{3}{2}, \\
\mathcal{O}_{g_0}\left(e^{-r}\right) & \text{if} \quad
a > \frac{3}{2}.
\end{cases}
\end{equation}
It is uniquely characterised by the fact that $\eta^0(\xi_0) = 1$ and $\gamma(\xi_0,\xi_0) = 0$.
\end{proposition}

\begin{proof}
Define $\bar{g}_0 = \eta^0\otimes \eta^0 + \gamma$, which is a continuous Riemannian metric on $\dK$ according to Theorem \ref{theorem:metric_expansion}.
Consider the continuous line bundle $\bar{L} = (\ker \eta^0)^{\perp_{\bar{g}_0}}$ on $\dK$.
The restriction of $\eta^0$ trivialises $\bar{L}$, which thus has a continuous nowhere vanishing section $\xi$.
Define $\xi_0 = \frac{\xi}{\eta^0(\xi)}$, which is continuous by construction.
Let $\{\eta^0,\ldots,\eta^{2n}\}$ be the limit coframe associated with any admissible frame.
Then $\eta^0(\xi_0) = 1$ and $\eta^j(\xi_0) = 0$ for $j\in \{1,\ldots,2n\}$.
In particular, $\xi_0$ is uniquely characterised by the relations $\eta^0(\xi_0)=1$ and $\gamma(\xi_0,\xi_0)=\sum_{j=1}^{2n}\eta^j(\xi_0)^2 = 0$.
Notice that for $j\in \{1,\ldots,2n\}$ and $r \geqslant 0$, one has
\begin{equation}
\label{eq:trick_eta}
\eta^j_r(\xi_0 - \xi_0^r) = \eta^j_r(\xi_0^r) - \eta^j_r(\xi) = \delta^j_0 - \eta^j_r(\xi_0) = \eta^j(\xi_0) - \eta^j_r(\xi_0)=  (\eta^j-\eta^j_r)(\xi_0),
\end{equation}
where $\delta$ stands for the Kronecker symbol.
Corollary \ref{corollary:pinching_g_g_r} yields the existence of a constant $c > 0$ such that $\|\xi_0^r - \xi_0\|_{g_0} \leqslant c e^{-\frac{r}{2}}\|Y_{(\xi_0^r - \xi_0)}\|_g$ for all $r \geqslant 0$.
The triangle inequality together with equation \eqref{eq:trick_eta} now yield
\begin{equation}
\|Y_{(\xi_0^r - \xi_0)}\|_g \leqslant
\big(e^r \|\eta^0-\eta^0_r\|_{g_0} + e^{\frac{r}{2}}\sum_{j=1}^{2n}\|\eta^j-\eta^j_r\|_{g_0}\big) \|\xi_0\|_{g_0}.
\end{equation}
Estimates \eqref{eq:estimates_xi_0} now follow from the estimates of Proposition \ref{proposition:eta_r_convergence}, together with the fact that $\|\xi_0\|_{g_0}$ is uniformly bounded by continuity of $\xi_0$ and compactness of $\dK$.
\end{proof}

\begin{remark}
\label{remark:xi_j}
Fix an admissible frame on $U\subset \dK$.
If $\xi_j^r = \mathcal{E}_r^* (e^{\frac{r}{2}}E_j)$ and if $\{\xi_0,\ldots,\xi_{2n}\}$ is the dual frame of $\{\eta^0,\ldots,\eta^{2n}\}$, a similar study shows that
\begin{equation}
\forall j \in \{1,\ldots,2n\},\quad  \xi_j - \xi_j^r =
\begin{cases}
\mathcal{O}_{g_0}\left(e^{-(a-\frac{1}{2})r}\right) & \text{if} \quad
1 < a <\frac{3}{2}, \\
\mathcal{O}_{g_0}\left((r+1)e^{-r}\right) & \text{if} \quad
a = \frac{3}{2}, \\
\mathcal{O}_{g_0}\left(e^{-r}\right) & \text{if} \quad
a > \frac{3}{2}.
\end{cases}
\end{equation}
The constants involved in the upper bounds are independent of the choice of the admissible frame and of $U$.
It relies on the fact that one can uniformly bound $\|\xi_j\|_{g_0}$ if $j\in \{1,\ldots,2n\}$, for instance, as an application of Corollary \ref{corollary:pinching_g_g_r}.
\end{remark}

For $v$ a vector field on $\dK$, the associated normal Jacobi fields $Y_v$ satisfies $\Ndr Y_v = SY_v$.
It follows from equation \eqref{eq:expression_jacobi_field} that in an admissible frame, one has
\begin{equation}
\label{eq:SY_v}
SY_v = \left(\delr \eta^0_r(v) + \eta^0_r(v) \right)e^r E_0
+ \sum_{j=1}^{2n}\left(\delr\eta^j_r(v) + \frac{1}{2}\eta^j_r(v) \right)e^{\frac{r}{2}}E_j.
\end{equation}
For $r \geqslant 0$, consider the pull-back $S_r = \mathcal{E}_r^*S$ of the shape operator $S$ through the diffeomorphism $\mathcal{E}_r \colon \dK \to \dK_r$.
It is well defined since $S$ leaves stable the tangent bundle of the level hypersurfaces $\dK_r$.

\begin{proposition}
\label{proposition:S_r_estimates}
Let $(M,g,J)$ be a complete, non-compact, almost Hermitian manifold of dimension at least $4$, admitting an essential subset $K$.
Assume that it satisfies the \ALCH and \AK conditions of order $a > \frac{1}{2}$.
Then the family $(S_r)_{r \geqslant 0}$ satisfies the estimates
\begin{equation}
\label{eq:Sr-truc}
S_r - \frac{1}{2}(\Id + \eta^0_r \otimes \xi_0^r) =
\begin{cases}
\mathcal{O}_{g_0}\left(e^{-(a-\frac{1}{2})r}\right) & \text{if} \quad
\frac{1}{2} < a <\frac{3}{2},\\
\mathcal{O}_{g_0}\left((r+1)e^{-r}\right) & \text{if} \quad
a = \frac{3}{2},\\
\mathcal{O}_{g_0}\left(e^{-r}\right) & \text{if} \quad
a > \frac{3}{2},
\end{cases}
\end{equation}
In particular, if $a > 1$, then $S_r \underset{r \to \infty}{\longrightarrow} \frac{1}{2}(\Id + \eta^0 \otimes \xi_0)$, and one can substitute $\eta^0_r\otimes \xi_0^r$ with $\eta^0 \otimes \xi_0$ in estimates \eqref{eq:Sr-truc}.
\end{proposition}

\begin{proof}
Let $v$ be a vector field on $\dK$.
It follows from Proposition \ref{proposition:eta_r_convergence} and from Corollary \ref{corollary:Y_v_asymptotic} that
\begin{equation}
SY_v -\frac{1}{2}(Y_v + \eta^0_r(v)e^rE_0) = \begin{cases}
\mathcal{O}_g\left(\|v\|_{g_0} e^{-(a-1)r}\right) & \text{if} \quad
\frac{1}{2} < a <\frac{3}{2},\\
\mathcal{O}_g\left(\|v\|_{g_0} (r+1)e^{-\frac{r}{2}}\right) & \text{if} \quad
a = \frac{3}{2},\\
\mathcal{O}_g\left(\|v\|_{g_0} e^{-\frac{r}{2}}\right) & \text{if} \quad
a > \frac{3}{2},\\
\end{cases}
\end{equation}
By the very definition of $S_r$, $\xi_0^r$ and $g_r$, it follows that
\begin{equation}
\big\|S_r-\frac{1}{2}(\Id + \eta^0_r\otimes \xi_0^r)\big\|_{g_r} =
\begin{cases}
\mathcal{O}\left(e^{-(a-1)r}\right) & \text{if} \quad
\frac{1}{2} < a <\frac{3}{2},\\
\mathcal{O}\left((r+1)e^{-\frac{r}{2}}\right) & \text{if} \quad
a = \frac{3}{2},\\
\mathcal{O}\left(e^{-\frac{r}{2}}\right) & \text{if} \quad
a > \frac{3}{2},\\
\end{cases}
\end{equation}
Finally, Corollary \ref{corollary:pinching_g_g_r} implies that
\begin{equation}
S_r - \frac{1}{2}(\Id + \eta^0_r \otimes \xi_0^r)
= \mathcal{O}_{g_0}\left(e^{-\frac{r}{2}}\big\|S_r - \frac{1}{2}(\Id + \eta^0_r \otimes \xi_0^r) \|_{g_r}\right),
\end{equation}
and estimates \eqref{eq:Sr-truc} now follow.
If $a > 1$, then estimates on $\|\eta^0-\eta^0_r\|_{g_0}$ (Proposition \ref{proposition:eta^0_global})
and on $\|\xi_0-\xi_0^r\|_{g_0}$ (Proposition \ref{proposition:xi_0_estimates}), together with the triangle inequality, show that one can replace $\eta^0_r\otimes \xi_0^r$ with $\eta^0\otimes \xi_0$ in estimates \eqref{eq:Sr-truc}.
This concludes the proof.
\end{proof}

\begin{remark}
In the complex hyperbolic space, the shape operator of a geodesic sphere of radius $r$, with outward unit normal $\nu$, is given by $S = \cotanh(r)\Id_{\R J\nu} + \frac{1}{2}\cotanh(\frac{r}{2}) \Id_{\{\nu,J\nu\}^{\perp}}$.
Proposition \ref{proposition:S_r_estimates} implies that the local extrinsic geometry of the level hypersurfaces $\dK_r$ is then asymptotic to that of horospheres in the complex hyperbolic space.
\end{remark}


\section{The almost complex structure}
\label{section:almost_complex}

\noindent This section is dedicated to prove the existence of a natural almost complex structure $J_0$ on the distribution of hyperplanes $H_0 = \ker \eta^0$, obtained as the restriction of a naturally defined tensor $\varphi$ on $\dK$.

The ambient almost complex structure $J$ does not leave stable the ambient distribution of hyperplanes $\{\delr\}^{\perp}$.
Consider the orthogonal projection $\pi \colon T \bar{M\setminus K} \to T \bar{M\setminus K}$ onto $\{\delr\}^{\perp}$.
Define $\Phi$ to be the field of endomorphisms on $\bar{M\setminus K}$ defined by $\Phi = \pi J \pi$.
Since $\pi$ and $J$ have unit norms, then $\|\Phi\|_g \leqslant 1$.
Formally, one has $\pi = \Id - g(\delr,\cdot) \otimes \delr$, and $\Phi$ then reads $\Phi = J + g(\cdot,\Jdelr) \otimes \delr - g(\cdot,\delr)\otimes \Jdelr$.

\begin{lemma}
\label{lemma:properties_phi}
Assume that $(M,g,J)$ satisfies the \AK condition of order $a > 0$.
For any admissible frame $\{E_0,\ldots,E_{2n}\}$ and any vector fields $X$ and $Y$, one has:
\begin{enumerate}
\item $g(\Phi X,\Phi Y) = g(X,Y) - g(X,\delr)g(Y,\delr) - g(X,\Jdelr)g(Y,\Jdelr)$,
\item $\Phi(E_0) = \mathcal{O}_g(e^{-ar})$,
\item $\Phi(E_j) - \JE_j = \mathcal{O}_g(e^{-ar})$ if $j\in \{1,\ldots,2n\}$.
\end{enumerate}
\end{lemma}

\begin{proof}
The first point is a straightforward computation.
To prove the second point, note that $\Phi(\Jdelr) = 0$, so that $\|\Phi(E_0)\|_g = \|\Phi(E_0-\Jdelr)\|_g \leqslant \|E_0-\Jdelr\|_g$.
The result follows from Corollary \ref{corollary:frame_estimates}.
Finally, by the very definition of $\Phi$, $\Phi(E_j)=\JE_j - g(E_j,\Jdelr)$, and the last point follows from Corollary \ref{corollary:frame_estimates}.
\end{proof}

The tensor $\Phi$ leaves stable the tangent distribution $\{\delr,\Jdelr\}^{\perp}$.
Therefore, one can pull it back through the family of diffeomorphisms $\left(\mathcal{E}_r\right)_{r\geqslant 0}$.

\begin{definition}
The family of endomorphisms $(\varphi_r)_{r \geqslant 0}$ is defined by $\varphi_r = \mathcal{E}_r^*\Phi$ for $r \geqslant 0$.
\end{definition}

Recall that $(S_r)_{r \geqslant 0}$ is the family of endomorphisms $\mathcal{E}_r^*S$ induced by the shape operator.

\begin{lemma}
\label{lemma:phi_r_estimates}
Assume that $(M,g,J)$ satisfies the \ALCH and \AK assumption of order $a > 1$.
Then the  following estimates hold:
\begin{enumerate}
\item $\varphi_r\xi_0^r = \mathcal{O}_{g_0}\left(e^{-(a-\frac{1}{2})r}\right)$.

\item $\varphi_r = \mathcal{O}_{g_0}(1)$,

\item $\eta^0_r\circ \varphi_r = \mathcal{O}_{g_0}(e^{-ar})$,

\item $\gamma_r(\varphi_r\cdot,\varphi_r\cdot) - \gamma_r = \mathcal{O}_{g_0}(e^{-(a-1)r})$,

\item $\varphi_r S_r - S_r \varphi_r =
\begin{cases}
\mathcal{O}_{g_0}(e^{-(a-\frac{1}{2})r}) & \text{if} \quad
1 < a <\frac{3}{2},\\
\mathcal{O}_{g_0}((r+1)e^{-r}) & \text{if} \quad
a = \frac{3}{2},\\
\mathcal{O}_{g_0}(e^{-r}) & \text{if} \quad
a > \frac{3}{2}.
\end{cases}
$
\end{enumerate}
\end{lemma}

\begin{proof}
We first show the first point.
From Corollary \ref{corollary:pinching_g_g_r}, there exists $c > 0$ such that for $r \geqslant 0$, $\|\varphi_r\xi_0^r\|_{g_0} \leqslant c \|\Phi (e^rE_0)\|_g e^{-\frac{r}{2}} = c\|\Phi (E_0)\|_g e^{\frac{r}{2}}$.
The result now follows from Lemma \ref{lemma:properties_phi}

Let us now focus on the second point.
Let $v$ be a vector field on $\dK$.
Corollary \ref{corollary:pinching_g_g_r} states that there exists $c>0$ such that $\|\varphi_rv\|_{g_0} \leqslant c \|\Phi(Y_v)\|_g e^{-\frac{r}{2}}$,
for all $r \geqslant 0$.
The result follows from the fourth point of Lemma \ref{lemma:properties_phi}.

For the third point, let $v$ be a vector field on $\dK$.
In an admissible frame, one has $\Phi(Y_v) = \eta^0_r(v) e^r \Phi(E_0) + e^{\frac{r}{2}}\sum_{j=1}^{2n}\eta^j_r(v) \Phi(E_j)$.
It then follows that
\begin{equation}
(\eta^0_r\circ \varphi_r)(v) = \eta^0_r(v) g(\Phi(E_0),E_0) + e^{-\frac{r}{2}}\sum_{j=1}^{2n} \eta^j_r(v) g(\Phi(E_j), E_0).
\end{equation}
Notice that $\Phi$ has range in $\{\Jdelr\}^{\perp}$, so that $g(\Phi(E_j), E_0)) = g(\Phi(E_j), E_0-\Jdelr)$ for all $j\in\{0,\ldots,2n\}$.
Recall that $\|\Phi\|_g \leqslant 1$ and that $\|E_j\|_g=1$ for all $j\in \{0,\ldots,2n\}$.
The triangle inequality now yields
$
\|\eta^0_r\circ \varphi_r\|_{g_0} \leqslant (\|\eta^0_r\|_{g_0}
+ e^{-\frac{r}{2}}\sum_{j=1}^n \|\eta^j_r\|_{g_0}) \|E_0-\Jdelr\|_g
$
for all $r \geqslant 0$.
The result follows from Corollary \ref{corollary:frame_estimates} (estimates on $E_0-\Jdelr$) and from Corollary \ref{corollary:eta_r_uniform_bound} (uniform bounds on $\{\|\eta^j_r\|_{g_0}\}_{j \in \{0,\ldots,2n\}}$).

Let us now consider the fourth point.
Let $u$ and $v$ be vector fields on $\dK$, and fix $r \geqslant 0$.
By Lemma \ref{lemma:properties_phi}, one has
$g_r(\varphi_ru,\varphi_rv) = g(\Phi Y_u,\Phi Y_v) = g(Y_u,Y_v) - g(Y_u,\Jdelr)g(Y_v,\Jdelr)$.
Cauchy-Schwarz inequality now yields
\begin{equation}
g_r(\varphi_ru,\varphi_rv) = g_r(u,v) - e^{2r}\eta^0_r(u)\eta^0_r(v) + \mathcal{O}(\|Y_u\|_g\|Y_v\|_g\|E_0-\Jdelr\|_g).
\end{equation}
It follows from Corollaries \ref{corollary:frame_estimates} and \ref{corollary:Y_v_asymptotic}, and from the very definition of $\gamma_r$, that
\begin{equation}
g_r(\varphi_r\cdot,\varphi_r\cdot) = e^r\gamma_r + \mathcal{O}_{g_0}( e^{(2-a)r}).
\end{equation}
Therefore, $e^{2r}(\eta^0_r\circ \varphi_r)\otimes(\eta^0_r\circ \varphi_r) + e^r \gamma_r(\varphi_r\cdot,\varphi_r\cdot) = e^r \gamma_r + \mathcal{O}_{g_0}(e^{(2-a)r})$.
From the preceding point, one has $e^{2r}(\eta^0_r\circ \varphi_r)\otimes(\eta^0_r\circ \varphi_r) = \mathcal{O}_{g_0}(e^{(2-2a)r})$, from which is deduced that $\gamma_r(\varphi_r\cdot,\varphi_r\cdot) = \gamma_r  + \mathcal{O}_{g_0}(e^{-(a-1)r})$
This concludes the proof of the fourth point.

Finally, let us prove the last point.
Write $S_r = S_r - \frac{1}{2}(\Id + \eta^0_r \otimes \xi_0^r) + \frac{1}{2}(\Id + \eta^0_r \otimes \xi_0^r)$, for $r \geqslant 0$.
By the triangle inequality, one has
\begin{equation}
\begin{split}
\|\varphi_r S_r - S_r \varphi_r \|_{g_0} &\leqslant 2 \|\varphi_r\|_{g_0}\|S_r -  \frac{1}{2}(\Id + \eta^0_r \otimes \xi_0^r)\|_{g_0} \\
& \quad +\frac{1}{2}(\|\eta^0_r\|_{g_0} \|\varphi_r\xi_0^r\|_{g_0} + \|\eta^0_r\circ \varphi_r\|_{g_0}\|\xi_0^r\|_{g_0}).
\end{split}
\end{equation}
The result now follows from uniform bounds on $\|\eta^0_r\|_{g_0}$ and $\|\xi_0^r\|_{g_0}$ (by uniform convergence), the estimates on $S_r - \frac{1}{2}(\Id + \eta^0_r \otimes \xi_0^r)$ (Proposition \ref{proposition:S_r_estimates}),
and the estimates on $\varphi_r$, $\eta^0_r\circ \varphi_r$,
and $\varphi_r \xi_0^r$, given by the three first points.
\end{proof}

We are now able to prove that the family $(\varphi_r)_{r \geqslant 0}$ converges to a continuous field of endomorphisms, provided that $a > 1$.

\begin{proposition}
\label{proposition:existence_phi}
Let $(M,g,J)$ be a complete, non-compact, almost Hermitian manifold of dimension at least $4$, with essential subset $K$.
Assume that it satisfies the \ALCH and \AK conditions of order $a > 1$.
Then there exists a continuous field of endomorphisms $\varphi$ on $\dK$ such that
\begin{equation}
\label{eq:phi_r-phi_estimates}
\varphi_r - \varphi =
\begin{cases}
\mathcal{O}_{g_0}\left(e^{-(a-\frac{1}{2})r}\right) & \text{if} \quad
1 < a <\frac{3}{2},\\
\mathcal{O}_{g_0}\left((r+1)e^{-r}\right) & \text{if} \quad
a = \frac{3}{2},\\
\mathcal{O}_{g_0}\left(e^{-r}\right) & \text{if} \quad
a > \frac{3}{2}.
\end{cases}
\end{equation}
In addition, $\varphi$ satisfies:
\begin{enumerate}
\item $\eta^0\circ \varphi = 0$ and $\varphi\xi_0 = 0$,
\item $\gamma(\varphi\cdot,\varphi\cdot) = \gamma$,
\item $\varphi^2 = -\Id + \eta^0 \otimes \xi_0$ and $\varphi^3 = -\varphi$.
\end{enumerate}
\end{proposition}

\begin{proof}
Let us first show the existence of $\varphi$.
The proof goes in two steps.
We first derive a differential equation for $(\varphi_r)_{r \geqslant 0}$.
Let $X$ be a vector field on $\bar{M\setminus K}$.
Then
\begin{equation}
\begin{split}
(\Ldr J)X &= [\delr,JX] - J[\delr,X]\\
&= (\Ndr(JX) - \nabla_{JX}\delr) - J(\Ndr X - \nabla_{X}\delr)\\
&= (\Ndr J) X + J \Ndr X - S(JX) - J \Ndr X + J(SX)\\
&= JSX - SJX + (\Ndr J)X.
\end{split}
\end{equation}
It follows that $\Ldr J = JS - SJ + \Ndr J$.
Recall that $\Phi = \pi J \pi$, where $\pi = \Id - g(\delr,\cdot)\otimes \delr$ is the orthogonal projection onto $\{\delr\}^{\perp}$.
It is a standard fact that $\Ldr g = 2g(S\cdot,\cdot)$.
Moreover, $S\delr = \Ndr \delr= 0$.
It follows that $\Ldr \pi = 0$, and consequently, that $\Ldr \Phi = \pi (JS - SJ + \Ndr J) \pi$.
Note that the eigenspaces of the projector $\pi$ are $\ker \pi = \R \delr$ and $\ker(\pi - \Id) = \{\delr\}^{\perp}$, which are both left stable by the shape operator $S$.
Hence, $S$ commutes with $\pi$, from which is derived that that $\Ldr \Phi = \Phi  S - S  \Phi + \pi (\Ndr J) \pi$.
Define $\psi_r = \mathcal{E}_r^*(\pi (\Ndr J) \pi)$, so that one has $\delr \varphi_r = \varphi_r S_r - S_r \varphi_r + \psi_r$.
A direct application of the \AK assumption and Corollary \ref{corollary:pinching_g_g_r} yields $\psi_r= \mathcal{O}_{g_0}(e^{-(a-\frac{1}{2})r})$.
Therefore, it follows from Lemma \ref{lemma:phi_r_estimates} that
\begin{equation}
\delr \varphi_r =
\begin{cases}
\mathcal{O}_{g_0}\left(e^{-(a-\frac{1}{2})r}\right) & \text{if} \quad
\frac{1}{2} < a <\frac{3}{2},\\
\mathcal{O}_{g_0}\left((r+1)e^{-r}\right) & \text{if} \quad
a = \frac{3}{2},\\
\mathcal{O}_{g_0}\left(e^{-r}\right) & \text{if} \quad
a > \frac{3}{2}.
\end{cases}
\end{equation}
Consequently, $(\varphi_r)_{r \geqslant 0}$ uniformly converges to some continuous tensor $\varphi$, which satisfies the inequality $\|\varphi_r - \varphi\|_{g_0} = \|\int_r^{\infty} \delr \varphi_r\|_{g_0} \leqslant \int_r^{\infty}\|\delr \varphi_r\|_{g_0}$ for all $r \geqslant 0$.
This implies estimates \eqref{eq:phi_r-phi_estimates}.

Let us now establish the claimed properties satisfied by $\varphi$.
The first two points are immediate consequences of Lemma \ref{lemma:phi_r_estimates}.
We thus focus on the last claim.
One easily checks that $\Phi$ satisfies the equality
$\Phi^2 = -\Id + g(\cdot,\Jdelr) \otimes \Jdelr + g(\cdot,\delr) \otimes \delr$.
Hence, one has ${\varphi_r}^2 = - \Id + \eta^0_r \otimes \xi_0^r + \epsilon_r$, for all $r \geqslant 0$, where $\epsilon_r = \mathcal{E}_r^*(g(\cdot,\Jdelr - E_0) \otimes \Jdelr + g(\cdot,E_0)\otimes (\Jdelr - E_0))$.
As usual, Corollary \ref{corollary:pinching_g_g_r} yields that
$\|\epsilon_r\|_{g_0} = \mathcal{O}(e^{\frac{r}{2}}\|E_0-\Jdelr\|_g) = \mathcal{O}(e^{-(a-\frac{1}{2})r})$, where the last equality is due to Corollary \ref{corollary:frame_estimates}.
The first part of the result now follows from the convergence of $(\eta^0_r)_{r \geqslant 0}$ and of $(\xi_0^r)_{r\geqslant 0}$ when $a > 1$.
The second part of the claim is a consequence of the first point.
\end{proof}

Proposition \ref{proposition:existence_phi} implies that when $a > 1$, $(\dK,\eta^0,\varphi,\xi_0)$ is an \emph{almost contact} manifold (see \cite{blair_riemannian_2010} for an introduction to this notion).
In particular, $\varphi$ induces an almost complex structure on the distribution of hyperplanes $H_0 = \ker \eta^0$.
The study conducted in this section finally implies the second of our main Theorems.

\begin{thm}
\label{theorem:almost_complex}
Let $(M,g,J)$ be a complete, non-compact almost Hermitian manifold of dimension greater than or equal to $4$
Assume that $M$ satisfies the \ALCH and \AK conditions of order $a > 1$.
Let $\eta^0$ and $\gamma$ be given by Theorem \ref{theorem:metric_expansion}, and let $\varphi$ be defined as in Proposition \ref{proposition:existence_phi}.
The restriction $J_0= \varphi|_{H_0}$ of $\varphi$ to the hyperplane distribution $H_0 = \ker \eta^0$ then induces an almost complex structure, and $\gamma^0=\gamma|_{H_0\times H_0}$ is $J_0$-invariant.
\end{thm}


\section{Higher regularity}
\label{section:higher_regularity}

\noindent This section is dedicated to show that under the stronger conditions \ALCHplus and \AKplus of order $a>1$, the tensors $\eta^0$, $\gamma$, and $\varphi$ defined previously gain in regularity.
As a consequence, we highlight a strictly pseudoconvex CR structure related to the expansion of the metric near infinity.

\subsection{Order one estimates}
We first provide several estimates that will be useful in the following study.

\begin{lemma}
\label{lemma:nabla_Yu_v}
Let $(M,g,J)$ be a complete, non-compact, almost Hermitian manifold of dimension at least $4$, admitting an essential subset $K$.
Assume that it satisfies the \ALCH condition of order $a > \frac{1}{2}$.
Let $u$ and $v$ be vector fields on $\dK$.
Let $V$ be the parallel transport of $v$ along radial geodesics.
Then $\nabla_{Y_u} V = \mathcal{O}_g(\|u\|_{g_0} \|v\|_{g_0} e^r)$.
\end{lemma}

\begin{proof}
Since $\Ndr V = 0$ and $[\delr,Y_u]=0$, one has $\Ndr (\nabla_{Y_u}V) = -R(\delr,Y_u)V$.
Hence, Kato's inequality yields $\big| \delr\|\nabla_{Y_u}V\|_g \big| \leqslant \|R\|_g \|Y_u\|_g \|V\|_g$ almost everywhere.
Recall that $\|R\|_g= \mathcal{O}(1)$ (Remark \ref{remark:R_bounded}) and that $\|V\|_g = \|v\|_{g_0}$.
Under the \ALCH condition of order $a > \frac{1}{2}$, one has $\|Y_u\|_g = \mathcal{O}(\|u\|_{g_0}e^r)$ (Corollary \ref{corollary:Y_v_asymptotic}).
The result follows from a straightforward integration.
\end{proof}

\begin{lemma}
\label{lemma:nala_Y_u_jdelr}
Let $(M,g,J)$ be a complete, non-compact, almost Hermitian manifold of dimension at least $4$, admitting an essential subset $K$.
Assume that it satisfies the \ALCH and \AK conditions of order $a > \frac{1}{2}$.
Then $\nabla_{Y_u}\Jdelr = \mathcal{O}_g(\|u\|_{g_0}e^r)$.
\end{lemma}

\begin{proof}
Write $\nabla_{Y_u}\Jdelr = (\nabla_{Y_u}J)\delr + J SY_u$.
Then $\|\nabla_{Y_u}\Jdelr\|_g \leqslant (\|\nabla J\|_g+ \|S\|_g) \|Y_u\|_g$, and the result follows from Lemma \ref{lemma:S_bounded}, the \AK assumption and the estimates of Corollary \ref{corollary:Y_v_asymptotic}.
\end{proof}

\begin{lemma}
\label{lemma:Ndr_nabla_J}
Assume that $(M,g,J)$ satisfies the \ALCH and \AKplus conditions of order $a > \frac{1}{2}$.
Then $\nabla_{Y_u}(\Ndr \Jdelr) = \mathcal{O}_g(\|u\|_{g_0}e^{-(a-1)r})$.
\end{lemma}

\begin{proof}
Since $\Ndr \delr = 0$ and $\nabla_{Y_u}\delr = SY_u$, it follows that
\begin{equation}
\begin{split}
\nabla_{Y_u}(\Ndr (J\delr)) & = \nabla_{Y_u}\left((\Ndr J)\delr\right) \\
&= \left(\nabla_{Y_u}(\Ndr J)\right)\delr + (\Ndr J) \nabla_{Y_u}\delr \\
&= (\nabla^2_{Y_u,\delr}J)\delr + (\nabla_{\nabla_{Y_u}\delr} J)\delr + (\Ndr J)\nabla_{Y_u}\delr\\
&= (\nabla^2_{Y_u,\delr}J)\delr + (\nabla_{SY_u}J)\delr + (\Ndr J)SY_u.
\end{split}
\end{equation}
The result follows from Corollary \ref{corollary:Y_v_asymptotic} (estimates on $SY_u$) and from the \AKplus assumption.
\end{proof}


\begin{lemma}
\label{lemma:nabla_Yu_v_estimates}
Assume that $(M,g,J)$ satisfies the \ALCHplus and \AK conditions of order $a > \frac{1}{2}$.
Let $\pi$ be the orthogonal projection onto $\{\delr\}^{\perp}$.
For $u$ and $v$ vector fields on $\dK$, one has:
\begin{enumerate}
\item $\pi((\nabla_{Y_u}S)Y_v) = \mathcal{O}_g(\|u\|_{g_0}\|v\|_{g_0}e^{\frac{3}{2}r})$.
\item $\pi(\nabla_{Y_u}Y_v) = \mathcal{O}_g\left((\|v\|_{g_0}+\|\nabla^{g_0}v\|_{g_0})\|u\|_{g_0}e^{\frac{3}{2}r} \right)$.
\end{enumerate}
\end{lemma}

\begin{proof}
We first consider the first point.
By Kato's inequality, and noticing that $\Ndr \pi = 0$, one has the almost everywhere inequality $\delr\|\pi(\nabla_{Y_u}S)Y_v)\|_g \leqslant \|\pi(\Ndr ((\nabla_{Y_u}S)Y_u))\|_g$.
The shape operator $S$ satisfies the Riccati equation $\Ndr S = -S^2 - R(\delr,\cdot)\delr$.
Moreover, one has $\pi S = S \pi$.
Direct computations using the equalities $\Ndr Y_v = SY_v$ and $\Ndr(SY_v) = -R(\delr,Y_v)\delr$ now yield
\begin{equation}
\begin{split}
\Ndr(\pi ((\nabla_{Y_u}S)Y_v))) & = \pi SR(\delr,Y_u)Y_v - \pi R(\delr,Y_u)SY_v - \pi R(SY_u,Y_v)\delr \\
& \quad -\pi R(\delr,Y_v)SY_u - \pi (\nabla_{Y_u}R)(\delr,Y_v)\delr - S \pi (\nabla_{Y_u}S)Y_v \\
& = \mathfrak{R} - S(\pi ((\nabla_{Y_u}S)Y_v))),
\end{split}
\end{equation}
where $\mathfrak{R}$ contains all the curvature terms.
From this is deduced the almost everywhere inequality $\delr (e^{-r}\|\pi ((\nabla_{Y_u}S)Y_v))\|_g) \leqslant e^{-r}\|\mathfrak{R}\|_g + (\|S\|_g-1) (e^{-r}\|\pi ((\nabla_{Y_u}S)Y_v))\|_g)$.
After a straightforward integration, Grönwall's Lemma yields
\begin{equation}
e^{-r}\|\pi ((\nabla_{Y_u}S)Y_v))\|_g \leqslant \left(\|(\nabla^{g}_uS)v\|_g + \int_{0}^r e^{-s}\|\mathfrak{R}\|_g \dx s\right)\exp\left(\int_0^r (\|S\|_g-1) \dx s\right).
\end{equation}
By tensoriality and compactness of $\dK$, one has $\|(\nabla^{g}_uS)v\|_g = \mathcal{O}(\|u\|_{g_0}\|v\|_{g_0})$.
Moreover, Lemma \ref{lemma:S_bounded} yields the estimate $\exp\left(\int_0^r (\|S\|_g-1)\dx s\right) = \mathcal{O}(1)$.
To conclude, it suffices to show that $\mathfrak{R} = \mathcal{O}_g(\|u\|_{g_0}\|v\|_{g_0}e^{\frac{3}{2}r})$.
The \ALCHplus assumption of order $a > \frac{1}{2}$ yields
\begin{equation}
\begin{split}
\mathfrak{R} &= \pi SR^0(\delr,Y_u)Y_v - \pi R^0(\delr,Y_u)SY_v - \pi R^0(SY_u,Y_v)\delr \\
& \quad -\pi R^0(\delr,Y_v)SY_u + \mathcal{O}_g\left( \|u\|_{g_0}\|v\|_{g_0}e^{-(a-2)r}\right).
\end{split}
\end{equation}
A close look at the definition of $R^0$ (see equation \eqref{eq:def_R0}) shows that the leading terms in $\|\mathfrak{R}\|_g$ are of the form $c\eta^0(u)\eta^j(v)e^{\frac{3}{2}r}$ or $c\eta^0(v)\eta^j(u)e^{\frac{3}{2}r}$ for $c$ a constant and $j \in \{1,\ldots,2n\}$.
The result follows.

Let us now show the second point.
Similarly, Kato's inequality yields the almost everywhere inequality
$\delr \|\pi(\nabla_{Y_u}Y_v)\|_g \leqslant \|\Ndr(\pi(\nabla_{Y_u}Y_v))\|_g$.
Straightforward computations, using that $\Ndr \pi = 0$, that $\pi$ and $S$ commute, and that $\Ndr Y_v = SY_v$, now yield the equality $\Ndr (\pi(\nabla_{Y_u}Y_v)) = -\pi R(Y_u,Y_v)\delr + \pi ((\nabla_{Y_u}S)Y_u) + S \pi (\nabla_{Y_u}Y_v)$.
Hence, one has
\begin{equation}
\begin{split}
\delr (e^{-r}\|\pi(\nabla_{Y_u}Y_v)\|_g) & \leqslant e^{-r}\|\pi R(Y_u,Y_v)\delr\|_g + e^{-r}\|\pi((\nabla_{Y_u}S)Y_v)\|_g \\
& \quad + (\|S\|_g-1) (e^{-r}\|\pi(\nabla_{Y_u}Y_v)\|_g) \quad a.e.
\end{split}
\end{equation}
The rest of the proof goes similarly to that of the first point, using the estimates derived on $\|\pi((\nabla_{Y_u}S)Y_v)\|_g$.
The main difference is that the initial data here is not tensorial in $v$, but instead is $\|\pi (\nabla_uv)\|_g = \|\nabla^{g_0}_uv\|_{g_0} \leqslant \|\nabla^{g_0}v\|_{g_0}\|u\|_{g_0}$.
\end{proof}

\begin{remark}
If one considers the whole vector field $\nabla_{Y_u}Y_v$ instead, then one only has the estimates $\|\nabla_{Y_u}Y_v\|_g = \mathcal{O}((\|v\|_{g_0}+\|\nabla^{g}v\|_{g})\|u\|_{g_0}e^{2r})$.
Indeed, the radial component is given by $g(\nabla_{Y_u}Y_v,\delr) = -g(SY_u,Y_v) \simeq -\eta^0(u)\eta^0(v)e^{2r}$ when $\eta^0(u)$ and $\eta^0(v)$ do not vanish.
\end{remark}

\subsection{Regularity of the admissible frames}
We shall now show that under the \ALCH and \AKplus conditions of order $a > 1$, the vector field $e_0$, defined in Definition \ref{definition:admissible}, is actually of class $\mathcal{C}^1$.

\begin{proposition}
\label{proposition:e_0_smooth}
Let $(M,g,J)$ be a complete, non-compact, almost Hermitian manifold of dimension at least $4$, admitting an essential subset $K$.
Assume that it satisfies the \ALCH and \AKplus conditions of order $a > 1$.
Then the vector field $e_0$ is of class $\mathcal{C}^1$; admissible frames can be chosen to have the same regularity.
\end{proposition}

\begin{proof}
It suffices to show that the $1$-form $\beta$ defined in Section \ref{subsection:admissible_frame} is of class $\mathcal{C}^1$.
To do so, we shall show that $\beta(v)$ is a $\mathcal{C}^1$ function for any $\mathcal{C}^1$ vector field $v$.
We prove this later fact by showing that $(u(\beta_r(v)))_{r\geqslant 0}$ uniformly converges for any $\mathcal{C}^1$ vector fields $u$ and $v$ on $\dK$.
Let $u$ and $v$ be such vector fields, and $r \geqslant 0$.
Then $u(\beta_r(v)) = Y_u(g(\Jdelr,V)) = \nabla_{Y_u}(g(\Jdelr,V))$, where $V$ is the parallel transport of $v$ along radial geodesics.
Since $[\delr,Y_u] = 0$ and $\Ndr V = 0$, one has
\begin{equation}
\delr\left(u (\beta_r(v))\right) = \Ndr(\nabla_{Y_u}(g(\Jdelr,V))) = \nabla_{Y_u}(\Ndr(g(\Jdelr,V))),
\end{equation}
so that $\delr(u (\beta_r(v))) = g(\nabla_{Y_u}(\Ndr(\Jdelr)),V) + g(\Ndr(\Jdelr),\nabla_{Y_u}V)$.
It now follows that one has $|\delr(u (\beta_r(v)))| \leqslant \|\nabla_{Y_u}V\|_g\|\Ndr(\Jdelr)\|_g + \|V\|_g\|\nabla_{Y_u}(\Ndr(\Jdelr))\|_g$.
Recall that $\|S\|_g = \mathcal{O}(1)$ (Lemma \ref{lemma:S_bounded}), $\|V\|_g = \|v\|_{g_0}$, and $\|Y_u\|_g = \mathcal{O}(\|u\|_{g_0}e^r)$ (Corollary \ref{corollary:Y_v_asymptotic}).
It now follows from Lemma \ref{lemma:nabla_Yu_v}, Lemma \ref{lemma:Ndr_nabla_J}, and the \AK assumption, that
\begin{equation}
\label{eq:delr_u_beta}
\delr\left(u (\beta_r(v))\right)  = \mathcal{O}\left(\|u\|_{g_0}\|v\|_{g_0}e^{-(a-1)r} \right).
\end{equation}
Consequently, $\delr(u (\beta_r(v)))$ uniformly converges for any vector fields $u$ and $v$.
This concludes the proof.
\end{proof}

It what follows, we will need to differentiate expressions involving $\nabla_{Y_u}E_j$ in the radial direction, with $Y_u$ a normal Jacobi field and $E_j$ an element of an admissible frame.
At a first glance, this is \emph{a priori} justified only if $E_j$ is of class $\mathcal{C}^2$.
One could prove such regularity by requiring the stronger condition $\|\nabla^3 J\|_g = \mathcal{O}(e^{-ar})$.
It turns out that one needs not assume this last condition, as a consequence of the fact that $E_j$ is solution to the first order linear differential equation $\Ndr E_j=0$.
Indeed, let $\{r,x^1,\ldots,x^{2n+1}\}$ be Fermi coordinates\footnote{That is, $\{x^1,\ldots,x^{2n+1}\}$ are coordinates on $\dK$, and that if $(x^1,\ldots,x^{2n+1})$ corresponds to $p\in \dK$, then $(r,x^1,\ldots,x^{2n+1})$ corresponds to $\mathcal{E}(r,p)\in M$.}, and write $E_j = \sum_{i=1}^{2n+1}E_j^i \partial_i$.
Then $\{E_j^i\}$ are solutions to the ODE $(E^i_j)' + \sum_{k=1}^{2n+1}E_j^kS_k^i = 0$, with $(S_k^i)$ the components of the shape operator $S$.
As a consequence, one can consider elements of the form $\Ndr (\nabla_{Y_u} E_j)$ even though $E_j$ is only of class $\mathcal{C}^1$.
In fact, one has $\Ndr(\nabla_{Y_u} E_j) = -R(\delr,Y_u)E_j$.

\begin{corollary}
\label{corollary:nabla(E_0-jdelr)}
Let $(M,g,J)$ be a complete, non-compact, almost Hermitian manifold of dimension at least $4$, admitting an essential subset $K$.
Assume that it satisfies the \ALCH and \AKplus conditions of order $a > 1$.
Let $u$ be a vector field on $\dK$.
Then
\begin{equation}
\nabla_{Y_u}(E_0 - \Jdelr) = \mathcal{O}_g(\|u\|_{g_0}e^{-(a-1)r}).
\end{equation}
\end{corollary}

\begin{proof}
Let $u$ be a vector field on $\dK$, and $\{E_0,\ldots,E_{2n}\}$ be an admissible frame of class $\mathcal{C}^1$.
Equation \eqref{eq:E_0-Jdelr} yields that
$\nabla_{Y_u}(E_0-\Jdelr) = -\sum_{j=0}^{2n} u(\beta_r(e_j)) E_j + \sum_{j=0}^{2n} (\delta_{0j} - \beta_r(e_j)) \nabla_{Y_u}E_j$.
During the proof of Proposition \ref{proposition:e_0_smooth}, we have shown that $(\beta_r)_{r \geqslant 0}$ converges in $\mathcal{C}^1$ topology.
Hence,
\begin{equation}
\forall j \in \{0,\ldots,2n\},\quad   \lim_{r \to \infty} u (\beta_r(e_j)) = u \left( \lim_{r \to \infty} \beta_r(e_j)\right)  = u(\beta(e_j)) = u(\delta_{0j}) = 0.
\end{equation}
Therefore, $|u(\beta_r(e_j))| = |\int_r^{\infty} \delr (u(\beta_r(e_j)))| \leqslant \int_r^{\infty} |\delr (u(\beta_r(e_j)))|$ for $j \in \{0,\ldots,2n\}$ and $r \geqslant 0$.
It follows from equation \eqref{eq:delr_u_beta} that $u(\beta_r(e_j)) = \mathcal{O}(\|u\|_{g_0}e^{-(a-1)r})$.
Moreover, by Corollary \ref{corollary:frame_estimates}, one has $|\delta_{0j}-\beta_r(e_j)| = \mathcal{O}(e^{-ar})$.
Finally, Lemma \ref{lemma:nabla_Yu_v} yields $\nabla_{Y_u}E_j = \mathcal{O}_g(\|u\|_ge^r)$.
The result now follows.
\end{proof}

\subsection{The contact form and the Carnot metric}
We shall now show that if the \ALCHplus and \AKplus conditions of order $a>1$ are satisfied, then $\eta^0$ and $\gamma|_{H_0\times H_0}$ are of class $\mathcal{C}^1$ and that $\dx\eta^0(\cdot,\varphi\cdot) = \gamma$.
In particular, $\eta^0$ is contact.
These results are analogous to \cite[Theorems B \& C]{pinoy_asymptotic2021}, although we give slightly different and considerably shorter proofs here.
The main difference is that we prove the $\mathcal{C}^1$ convergence of elements of the form $(\eta^j_r(v))_{r\geqslant 0}$, instead of $\mathcal{C}^0$ convergence of elements of the form $(\mathcal{L}_u\eta^j_r)_{r\geqslant 0}$.

\begin{thm}
\label{theorem:contact}
Let $(M,g,J)$ be a complete, non-compact, almost Hermitian manifold of dimension at least $4$, with essential subset $K$.
Assume that it satisfies the \ALCHplus and \AKplus conditions of order $a > 1$.
Then $\eta^0$ is a contact form of class $\mathcal{C}^1$.
Moreover, $\dx \eta^0(\cdot,\varphi \cdot) = \gamma$, and the Reeb vector field of $\eta^0$ is $\xi_0$.
\end{thm}

\begin{proof}
The proof is divided in three parts.
First, we show that $\eta^0$ is of class $\mathcal{C}^1$.
Then we derive an expression for $\dx \eta^0(\cdot,\varphi\cdot)$, and deduce that $\eta^0$ is contact.
Finally, we show that $\xi_0$ is the Reeb vector field of $\eta^0$.

To show that $\eta^0$ is of class $\mathcal{C}^1$, we show that for any vector field $v$, the function $\eta^0(v)$ is of class $\mathcal{C}^1$.
To do so, we show that for any other vector field $u$, $(u(\eta^0_r(v)))_{r \geqslant 0}$ uniformly converges on $\dK$.
Let $u$ and $v$ be vector fields on $\dK$.
Let $f$ be the function on $\bar{M\setminus K}$ defined by the expression $f= e^r\left(u(\eta^0_r(v)\right) = Y_u\left(g(Y_v,E_0)\right) = \nabla_{Y_u}\left(g(Y_u,E_0) \right)$.
Then $f$ is smooth in the radial direction.
Since $[\delr,Y_u]=0$ and $\Ndr E_0=0$, one has
\begin{equation}
\delr f = \Ndr (\nabla_{Y_u} ((g(Y_v,E_0))) = \nabla_{Y_u}(\Ndr (g(Y_v,E_0))) = \nabla_{Y_u}(g(\Ndr Y_v,E_0)).
\end{equation}
Similarly, one has $\delr^2f = \nabla_{Y_u}(g(\Ndr(\Ndr Y_v),E_0))$.
For $Y_v$ is a Jacobi field, one has the equality $\Ndr(\Ndr Y_v) = -R(\delr,Y_v)\delr$, and it follows that $\delr^2f = -\nabla_{Y_u}(R(\delr,Y_v,\delr,E_0))$.
Notice that
\begin{equation}
\begin{split}
R(\delr,Y_v,\delr,E_0) &= R(\delr,Y_v,\delr,\Jdelr) + R(\delr,Y_v,\delr,E_0-\Jdelr)\\
&= R^0(\delr,Y_v,\delr,\Jdelr) + R(\delr,Y_v,\delr,E_0-\Jdelr)\\
& \quad + (R-R^0)(\delr,Y_v,\delr,\Jdelr).
\end{split}
\end{equation}
One readily checks from the definition of $R^0$ that $R^0(\delr,Y_v,\delr,\Jdelr) = -g(Y_v,\Jdelr)$, so that $R^0(\delr,Y_v,\delr,\Jdelr) = -g(Y_v,E_0) - g(Y_v,\Jdelr - E_0)$.
Hence, it follows that
\begin{equation}
\begin{split}
\delr^2f - f &= g(\nabla_{Y_u}Y_v,\Jdelr -E_0) + g(Y_v,\nabla_{Y_u}(\Jdelr-E_0)) \\
& \quad - (\nabla_{Y_u}R)(\delr,Y_v,\delr,E_0-\Jdelr) - R(SY_u,Y_v,\delr,E_0-\Jdelr) \\
& \quad - R(\delr,\nabla_{Y_u}Y_u,\delr,E_0-\Jdelr) - R(\delr,Y_v,SY_u,E_0-\Jdelr) \\
& \quad - R(\delr,Y_v,\delr,\nabla_{Y_u}(E_0-\Jdelr)) - (\nabla_{Y_u}(R-R^0))(\delr,Y_v,\delr,\Jdelr) \\
& \quad - (R-R^0)(SY_u,Y_v,\delr,\Jdelr) - (R-R^0)(\delr,\nabla_{Y_u}Y_v,\delr,\Jdelr) \\
& \quad - (R-R^0)(\delr,Y_v,SY_u,\Jdelr) - (R-R^0)(\delr,Y_v,\delr,\nabla_{Y_u}\Jdelr).
\end{split}
\end{equation}
Note that the radial part of $\nabla_{Y_u}Y_v$ plays no role here due to the symmetries of the Riemann curvature tensor, so that one can substitute $\nabla_{Y_u}Y_v$ with $\pi(\nabla_{Y_u}Y_v)$ in this latter expression.
Recall that one has the following estimates:
\begin{itemize}
\item $R, S = \mathcal{O}_g(1)$ (Remark \ref{remark:R_bounded} and Lemma \ref{lemma:S_bounded}),
\item $R-R^0,\nabla R, \nabla(R-R^0) = \mathcal{O}_g(e^{-ar})$ (\ALCHplus condition and Remark \ref{remark:nabla(R-R^0)}),
\item $E_0-\Jdelr = \mathcal{O}_g(e^{-ar})$ (Corollary \ref{corollary:frame_estimates}),
\item $Y_u,Y_v = \mathcal{O}_g(\|u\|_{g_0}e^r)$ (Corollary \ref{corollary:Y_v_asymptotic}),
\item $\nabla_{Y_u}\Jdelr = \mathcal{O}_g(\|u\|_{g_0}e^r)$ (Lemma \ref{lemma:nala_Y_u_jdelr}),
\item $\pi(\nabla_{Y_u}Y_v) = \mathcal{O}_g((\|v\|_{g_0}+\|\nabla^{g_0}v\|_{g_0})\|u\|_{g_0}e^{\frac{3}{2}r})$ (Lemma \ref{lemma:nabla_Yu_v_estimates}),
\item $\nabla_{Y_u}(E_0-\Jdelr) = \mathcal{O}_g(\|u\|_{g_0}e^{-(a-1)r})$ (Corollary \ref{corollary:nabla(E_0-jdelr)}).
\end{itemize}
Hence, the triangle inequality yields
\begin{equation}
\label{eq:y''-y}
\delr^2f - f = \mathcal{O}\left((\|v\|_{g_0}+\|\nabla^{g_0}v\|_{g_0})\|u\|_{g_0}e^{(2-a)r}\right).
\end{equation}
Define $h = \delr f - f$, and notice that $\delr h + h = \delr^2f - f$.
It now follows from equation \eqref{eq:y''-y} that $\delr (e^rh) = \mathcal{O}\left((\|v\|_{g_0}+\|\nabla^{g_0}v\|_{g_0})\|u\|_{g_0}e^{(3-a)r}\right)$.
Therefore, one has
\begin{equation}
e^rh =
\begin{cases}
\mathcal{O}\left((\|v\|_{g_0}+\|\nabla^{g_0}v\|_{g_0})\|u\|_{g_0}e^{(3-a)r}\right) & \text{if} \quad 1 < a < 3,\\
\mathcal{O}\left((\|v\|_{g_0}+\|\nabla^{g_0}v\|_{g_0})\|u\|_{g_0}(r+1)\right) & \text{if} \quad a=3,\\
\mathcal{O}\left((\|v\|_{g_0}+\|\nabla^{g_0}v\|_{g_0})\|u\|_{g_0}\right) & \text{if} \quad a > 3.
\end{cases}
\end{equation}
Notice that $e^{-r}h = \delr(e^{-r}f) = \delr \left(u\left(\eta^0_r(v)\right) \right)$.
Hence,
\begin{equation}
\delr \left(u\left(\eta^0_r(v)\right) \right) =
\begin{cases}
\mathcal{O}\left((\|v\|_{g_0}+\|\nabla^{g_0}v\|_{g_0})\|u\|_{g_0}e^{-(a-1)r}\right) & \text{if} \quad 1 < a < 3,\\
\mathcal{O}\left((\|v\|_{g_0}+\|\nabla^{g_0}v\|_{g_0})\|u\|_{g_0}(r+1)e^{-2r}\right) & \text{if} \quad a=3,\\
\mathcal{O}\left((\|v\|_{g_0}+\|\nabla^{g_0}v\|_{g_0})\|u\|_{g_0}e^{-2r}\right) & \text{if} \quad a > 3.
\end{cases}
\end{equation}
Consequently, $\left(u(\eta^0_r(v))\right)_{r\geqslant 0}$ uniformly converges as $r\to \infty$,
and $\eta^0$ is then of class $\mathcal{C}^1$.

We shall now derive an expression for $\dx \eta^0(\cdot,\varphi\cdot)$, by computing the limit of $\dx \eta^0_r(\cdot,\varphi_r\cdot)$ as $r \to \infty$.
Let $u$ and $v$ be vector fields on $\dK$.
For $r \geqslant 0$, it holds that
\begin{equation}
\begin{split}
\dx \eta^0_r(u,\varphi_rv) &= u\left(\eta^0_r(\varphi_rv)\right) - (\varphi_rv)\left(\eta^0_r(u)\right) - \eta^0_r([u,\varphi_rv])\\
&= e^{-r} \left( Y_u g(\Phi Y_v,E_0) - (\Phi Y_v)g(Y_u,E_0) - g([Y_u,\Phi Y_v],E_0) \right)\\
&= e^{-r}\left(g(\Phi Y_v,\nabla_{Y_u}E_0) - g(Y_u,\nabla_{\Phi Y_v}E_0)\right).
\end{split}
\end{equation}
On the one hand, it holds that
\begin{equation}
\begin{split}
g(\Phi Y_v,\nabla_{Y_u}E_0) &= g(\Phi Y_v,\nabla_{Y_u}\Jdelr) + g(\Phi Y_v,\nabla_{Y_u}(E_0-\Jdelr))\\
&= g(\Phi Y_v,JSY_u) + g(\Phi Y_v,(\nabla_{Y_u}J)\delr)+ g(\Phi Y_v,\nabla_{Y_u}(E_0-\Jdelr))\\
&= -g(J\Phi Y_v,SY_u) + g(\Phi Y_v,(\nabla_{Y_u}J)\delr)+ g(\Phi Y_v,\nabla_{Y_u}(E_0-\Jdelr)).
\end{split}
\end{equation}
On the other hand, one has
\begin{equation}
\begin{split}
g(Y_u,\nabla_{\Phi Y_v}E_0) &= g(Y_u,\nabla_{\Phi Y_v} \Jdelr) + g(Y_u,\nabla_{\Phi Y_v}(E_0-\Jdelr)) \\
&= g(Y_u,JS\Phi Y_v) + g(Y_u, (\nabla_{\Phi Y_v}J)\delr) + g(Y_u,\nabla_{\Phi Y_v}(E_0-\Jdelr))\\
&= -g(JY_u,S\Phi Y_v) + g(Y_u, (\nabla_{\Phi Y_v}J)\delr) + g(Y_u,\nabla_{\Phi Y_v}(E_0-\Jdelr)).
\end{split}
\end{equation}
It then follows from the \AK assumption, Corollary \ref{corollary:Y_v_asymptotic} and Corollary \ref{corollary:nabla(E_0-jdelr)} that
\begin{equation}
\label{eq:deta^0_r}
\dx\eta^0_r(u,\varphi_rv) = e^{-r}\left(g(JY_u,S\Phi Y_v) - g(J\Phi Y_v,SY_u)\right) + \mathcal{O}\left(\|u\|_{g_0}\|v\|_{g_0}e^{-(a-1)r}\right).
\end{equation}
Fix $\{E_0,\ldots,E_{2n}\}$ an admissible frame.
From Corollary \ref{corollary:frame_estimates} and Corollary \ref{corollary:Y_v_asymptotic}, one has the estimate $Y_v = \eta^0(v) e^r \Jdelr + \sum_{j=1}^{2n}\eta^j(v)e^{\frac{r}{2}}E_j + \mathcal{O}_g(\|v\|_{g_0}e^{-(a-1)r})$.
It now follows from Lemma \ref{lemma:properties_phi} that $J\Phi Y_v = -\sum_{j=1}^{2n} \eta^j(v) e^{\frac{r}{2}} E_j + \mathcal{O}_g(\|v\|_{g_0}e^{-(a-1)r})$.
Corollary \ref{corollary:Y_v_asymptotic} now yields
\begin{equation}
\label{eq:deta^0_part_1}
g(J\Phi Y_v,SY_u) = -\frac{e^r}{2} \sum_{j=1}^{2n} \eta^j(v)\eta^j(u) + \mathcal{O}(\|u\|_{g_0}\|v\|_{g_0}e^{-(a-2)r}).
\end{equation}
Similarly, one shows that
\begin{equation}
\label{eq:deta^0_part_2}
g(JY_u,S\Phi Y_v) = \frac{e^r}{2} \sum_{j=1}^{2n}\eta^j(u)\eta^j(v) + \mathcal{O}(\|u\|_{g_0}\|v\|_{g_0}e^{-(a-2)r}).
\end{equation}
Recall the local expression $\gamma = \sum_{j=1}^{2n}\eta^j\otimes \eta^j$.
Equations \eqref{eq:deta^0_r}, \eqref{eq:deta^0_part_1} and \eqref{eq:deta^0_part_2} now yield
\begin{equation}
\dx \eta^0_r(u,\varphi_rv) = \gamma(u,v) + \mathcal{O}(\|u\|_{g_0}\|v\|_{g_0}e^{-(a-1)r}).
\end{equation}
By uniform convergence of the first derivatives of $(\eta^0_r)_{r\geqslant 0}$, it follows that $\dx \eta^0(\cdot,\varphi \cdot) = \gamma$.
Proposition \ref{proposition:gamma_r_convergence} hence shows that $\dx \eta^0$ is non-degenerate on $\ker \eta^0$.
In particular, $\eta^0$ is a contact form.

To conclude, let us show that $\xi_0$ is the Reeb vector field of $\eta^0$.
Since $\eta^0(\xi_0) = 1$, it remains to show that $\dx \eta^0(\xi_0,v) = 0$ for all vector field $v$ tangent to $H_0$.
Let $v$ be such a vector field.
The image of $\varphi$ being exactly $H_0$, there exists a vector field $u$ on $\dK$ such that $v = \varphi u$.
By Proposition \ref{proposition:existence_phi}, $\gamma$ is $\varphi$-invariant and $\varphi \xi_0=0$.
From the preceding point, $\dx \eta^0(\cdot,\varphi\cdot) = \gamma$.
Hence, $\dx \eta^0(\xi_0,v) = \dx \eta^0(\xi_0,\varphi u) = \gamma(\xi_0,u) = \gamma(\varphi \xi_0,\varphi u) = \gamma(0,\varphi u) = 0$.
This concludes the proof.
\end{proof}

\begin{corollary}
Under the assumptions of Theorem \ref{theorem:contact}, the distribution $H_0 = \ker \eta^0$ is a contact distribution of class $\mathcal{C}^1$.
\end{corollary}

The next result shows that under the assumptions of Theorem \ref{theorem:contact}, the Carnot metric $\gamma^0$ on $H_0$ is of the same regularity.
The proof is very similar.

\begin{thm}
\label{theorem:gamma_c^1}
Let $(M,g,J)$ be a complete, non-compact, almost Hermitian manifold of dimension at least $4$, with essential subset $K$.
Assume that it satisfies the \ALCHplus and \AKplus conditions of order $a > 1$.
Then $\gamma^0 = \gamma|_{H_0\times H_0}$ is of class $\mathcal{C}^1$.
\end{thm}

\begin{proof}
Let $\{E_0,\ldots,E_{2n}\}$ be an admissible frame of class $\mathcal{C}^1$ defined on a cone $E(\R_+\times U)$, and fix $j\in \{1,\ldots,2n\}$.
Let us first show that $\eta^j$ is of class $\mathcal{C}^1$ on the distribution $H_0|_U$.
To do so, we shall prove that $\left(u\left(\eta^j_r(v)\right)\right)_{r \geqslant 0}$ locally uniformly converges on $U$ for $v$ tangent to $H_0|_U$ and $u$ any vector field on $U$.

Let $u$ and $v$ be such vector fields, and $r \geqslant 0$ be fixed.
Let $f^j = e^{\frac{r}{2}} \, u\left(\eta^j_r(v)\right) = \nabla_{Y_u}\left(g(Y_v,E_j)\right)$, which is smooth in the radial direction.
Since $[\delr,Y_u] = 0$ and $\Ndr E_j = 0$, one has
\begin{equation}
\delr^2 f^j = \Ndr(\Ndr(\nabla_{Y_u}\left(g(Y_v,E_j)\right))) = \nabla_{Y_u}\, g(\Ndr(\Ndr Y_v),E_j),
\end{equation}
and, for $Y_v$ is a Jacobi field, one has $\delr^2f^j = - \nabla_{Y_u}(R(\delr,Y_v,\delr,E_j))$.
One checks from the very definition of $R^0$ that $R^0(\delr,Y_v,\delr,E_j) = -\frac{1}{4}g(Y_v,E_j) - \frac{3}{4}g(Y_v,\Jdelr)g(E_j,\Jdelr)$.
Therefore, one has the equality
\begin{equation}
\begin{split}
\delr^2f^j - \frac{1}{4}f^j &= \frac{3}{4}g(\nabla_{Y_u}Y_v,\Jdelr)g(E_j,\Jdelr) + \frac{3}{4}g(Y_v,\nabla_{Y_u}\Jdelr)g(E_j,\Jdelr) \\
&\quad + \frac{3}{4}g(Y_v,\Jdelr)g(\nabla_{Y_u}E_j,\Jdelr) + \frac{3}{4}g(Y_v,\Jdelr)g(E_j,\nabla_{Y_u}\Jdelr)\\
& \quad - \nabla_{Y_u}(R-R^0)(\delr,Y_v,\delr,E_j) - (R-R^0)(SY_u,Y_v,\delr,E_j) \\
& \quad - (R-R^0)(\delr,\nabla_{Y_u}Y_v,\delr,E_j) - (R-R^0)(\delr,Y_v,SY_u,E_j)\\
& \quad - (R-R^0)(\delr,Y_v,\delr,\nabla_{Y_u}E_j).
\end{split}
\end{equation}
As in the proof of Theorem \ref{theorem:almost_complex}, the radial component of $\nabla_{Y_u}Y_v$ plays no role due to the symmetries of $R$, so that one can substitute this term with $\pi(\nabla_{Y_u}Y_v)$.
Moreover, $g(E_j,\Jdelr) = \beta_r(e_j)$, where $(\beta_r)_{r \geqslant 0}$ is the family defined in Section \ref{subsection:admissible_frame}.
Recall that one has the following estimates:
\begin{itemize}
\item $R, S = \mathcal{O}_g(1)$ (Remark \ref{remark:R_bounded} and Lemma \ref{lemma:S_bounded}),
\item $R-R^0,\nabla (R-R^0) = \mathcal{O}_g(e^{-ar})$, (\ALCHplus condition and Remark \ref{remark:nabla(R-R^0)}),
\item $\beta_r(e_j) = \mathcal{O}(e^{-ar})$ (Corollary \ref{corollary:frame_estimates}),
\item $Y_u = \mathcal{O}_g(\|u\|_{g_0}e^r)$ and $Y_v = \mathcal{O}_g(\|v\|_{g_0}e^\frac{r}{2})$ (Corollary \ref{corollary:Y_v_asymptotic}),
\item $\nabla_{Y_u}E_j = \mathcal{O}_g(\|u\|_{g_0}e^r)$ (Lemma \ref{lemma:nabla_Yu_v}),
\item $\nabla_{Y_u}\Jdelr = \mathcal{O}_g(\|u\|_{g_0}e^r)$ (Lemma  \ref{lemma:nala_Y_u_jdelr}),
\item $\pi(\nabla_{Y_u}Y_v) = \mathcal{O}_g((\|\nabla^{g_0}u\|_{g_0} + \|u\|_{g_0})\|v\|_{g_0}e^{\frac{3}{2}r})$
(Lemma \ref{lemma:nabla_Yu_v_estimates}).
\end{itemize}
It follows from the triangle inequality that $\delr^2 f^j - \frac{1}{4}f^j = \mathcal{O}((\|v\|_{g_0}+\|\nabla^{g_0}v\|_{g_0})\|u\|_{g_0}e^{-(a-\frac{3}{2})r})$.
Let $h^j$ be the function defined by $h^j = \delr f^j - \frac{1}{2}f^j$.
Then $\delr h^j + \frac{1}{2}h^j = \delr^2f^j - \frac{1}{4}f^j$, from which is derived that $\delr(e^{\frac{r}{2}}h^j) = \mathcal{O}((\|v\|_{g_0}+\|\nabla^{g_0}v\|_{g_0})\|u\|_{g_0}e^{-(a-2)r})$.
A straightforward integration now yields
\begin{equation}
e^{\frac{r}{2}}h^j =
\begin{cases}
\mathcal{O}\left((\|v\|_{g_0}+\|\nabla^{g_0}v\|_{g_0})\|u\|_{g_0}e^{(2-a)r}\right) & \text{if} \quad 1 < a < 2,\\
\mathcal{O}\left((\|v\|_{g_0}+\|\nabla^{g_0}v\|_{g_0})\|u\|_{g_0}(r+1)\right) & \text{if} \quad a = 2,\\
\mathcal{O}\left((\|v\|_{g_0}+\|\nabla^{g_0}v\|_{g_0})\|u\|_{g_0}\right) & \text{if} \quad a > 2.
\end{cases}
\end{equation}
Notice that $e^{-\frac{r}{2}}h^j = \delr(e^{-\frac{r}{2}}f^j) = \delr\left( u(\eta^j_r(v))\right)$, from which is deduced that
\begin{equation}
\delr\left( u (\eta^j_r(v))\right) =
\begin{cases}
\mathcal{O}\left((\|v\|_{g_0}+\|\nabla^{g_0}v\|_{g_0})\|u\|_{g_0}e^{-(a-1)r}\right) & \text{if} \quad 1 < a < 2,\\
\mathcal{O}\left((\|v\|_{g_0}+\|\nabla^{g_0}v\|_{g_0})\|u\|_{g_0}(r+1)e^{-r}\right) & \text{if} \quad a = 2,\\
\mathcal{O}\left((\|v\|_{g_0}+\|\nabla^{g_0}v\|_{g_0})\|u\|_{g_0}e^{-r}\right) & \text{if} \quad a > 2.
\end{cases}
\end{equation}
In any case, $\left( u(\eta^j_r(v))\right)_{r\geqslant 0}$ locally uniformly converges.
As a consequence, $\eta^j|_{H_0|_U}$ is of class $\mathcal{C}^1$.
We immediately deduce from the local expression $\gamma = \sum_{j=1}^{2n}\eta^j\otimes \eta^j$ that $\gamma^0=\gamma|_{H_0\times H_0}$ is of class $\mathcal{C}^1$.
This concludes the proof.
\end{proof}

\begin{remark}
\label{remark:a>3/2}
With the stronger assumption $a > \frac{3}{2}$, the same proof shows that for $j\in \{1,\ldots,2n\}$, $\eta^j$ is of class $\mathcal{C}^1$ in all directions, and so is $\gamma$.
Indeed, in this case, on has to consider the estimate $Y_v = \mathcal{O}_g(\|v\|_{g_0}e^r)$ instead.
\end{remark}

\subsection{The almost complex structure}
We shall now show that the almost complex structure $J_0$ defined on the $\mathcal{C}^1$ distribution $H_0$ is of the same regularity, and that it is formally integrable.
We first remark that the local vector fields $\{\xi_1,\ldots,\xi_{2n}\}$ are of class $\mathcal{C}^1$, although the Reeb vector field $\xi_0$ might only be continuous.

\begin{lemma}
\label{lemma:coframe_c1}
Let $(M,g,J)$ be a complete, non-compact, almost Hermitian manifold of dimension at least 4, with essential subset $K$.
Assume that $(M,g,J)$ satisfies the \ALCHplus and \AKplus conditions of order $a > 1$.
Let $\{\eta^0,\ldots,\eta^{2n}\}$ be the local coframe associated to any admissible frame $\{E_0,\ldots,E_{2n}\}$.
Let $\{\xi_0,\xi_1,\ldots,\xi_{2n}\}$ be its dual frame.
Then for $j\in \{1,\ldots,2n\}$, $\xi_j$ is a vector field of class $\mathcal{C}^1$.
\end{lemma}

\begin{proof}
Throughout the proof of Theorem \ref{theorem:contact}, we have shown that $\{\eta^1,\ldots,\eta^{2n}\}$ is a $\mathcal{C}^1$ trivialisation of the $\mathcal{C}^1$ vector bundle $\Hom(H_0,\R)$.
Consequently, $\{\xi_1,\ldots,\xi_{2n}\}$ is a $\mathcal{C}^1$ trivialisation of the vector bundle $H_0$.
\end{proof}

We now show that under the \AKplus condition of order $a > 0$, admissible frames can almost be chosen to be $J$-frames, in the following sense.

\begin{lemma}
\label{lemma:J_admissible_frame}
Let $(M,g,J)$ be a complete, non-compact, almost Hermitian manifold of dimension at least $4$, and with essential subset $K$.
Assume that it satisfies the \AKplus condition of order $a > 0$.
Then there exists an admissible frame $\{E_0,\ldots,E_{2n}\}$ such that
\begin{equation}
\forall j \in \{1,\ldots,n\},\quad \JE_{2j-1} - E_{2j} =  \mathcal{O}_g(e^{-ar}).
\end{equation}
\end{lemma}

\begin{proof}
Let $U\subset \dK$ be an open domain on which $H_0$ is trivialisable.
Let $e_1$ be a unit section of $H_0|_U$ of class $\mathcal{C}^1$, and let $E_1$ be its parallel transport along radial geodesics.
Consider the family of $1$-forms $\beta^1_r \colon {H_0}|_U \to \R$ defined by $\beta^1_r(v) = g(V, \JE_1)|_{\dK_r}$,
where $V$ is the parallel transport of $v$ along radial geodesics.
The same study than that conducted for the proofs of Lemma \ref{lemma:beta_r_converges} and Proposition \ref{proposition:e_0_smooth} shows that under the \AKplus condition of order $a >1$, there exists a nowhere vanishing $1$-form $\beta^1$ on $U$, which is of class $\mathcal{C}^1$, such that $\|\beta^1_r - \beta^1\|_{g_0} =\mathcal{O}(e^{-ar})$.
Let $e_2$ be the unique $\mathcal{C}^1$ section of $H_0|_U$ such that $e_2 \perp^{g_0} \ker \beta^1$, $\|e_2\|_{g_0} = 1$ and $\beta^1(e_2) > 0$.
Define $E_2$ to be its parallel transport along radial geodesics.
Similarly to Corollary \ref{corollary:frame_estimates}, one shows that $E_2-\JE_1 = \mathcal{O}_g(e^{-ar})$.
The rest of the proof follows by induction.
\end{proof}

We refer to such an admissible frame as a \emph{$J$-admissible frame}.
We are now able to show the last Theorem of this section, exhibiting a strictly pseudoconvex CR structure at infinity.

\begin{thm}
\label{theorem:strictly_pseudoconvex}
Let $(M,g,J)$ be a complete, non-compact, almost Hermitian manifold of dimension at last 4, with essential subset $K$.
Assume that it satisfies the \ALCHplus and \AKplus condition of order $a > 1$.
Let $J_0$ be the almost complex structure on $H_0$ induced by $\varphi$.
Then $J_0$ is of class $\mathcal{C}^1$, and is formally integrable.
In particular, $(\dK,H_0,J_0)$ is a strictly pseudoconvex CR manifold of class $\mathcal{C}^1$.
\end{thm}

\begin{proof}
Let $\{E_0,\ldots,E_{2n}\}$ be a $J$-admissible frame of class $\mathcal{C}^1$, and $\{\eta^1,\ldots,\eta^{2n}\}$ and $\{\xi_1,\ldots,\xi_{2n}\}$ be the associated $\mathcal{C}^1$ coframe and frame.
Then $\{\delr,E_0,\ldots,E_{2n}\}$ is an orthonormal frame.
Since $\Phi(\delr) = \Phi(\Jdelr)= 0$, one has $\Phi = \sum_{j=0}^{2n} g(\cdot,E_j)\otimes \Phi(E_j)$.
It then follows from Lemma \ref{lemma:properties_phi} and Lemma \ref{lemma:J_admissible_frame} that $\Phi = \sum_{j=1}^n g(\cdot,E_{2j-1})\otimes E_{2j} - g(\cdot,E_{2j})\otimes E_{2j-1} + \mathcal{O}_g(e^{-ar})$.
Corollary \ref{corollary:pinching_g_g_r} now yields
$\varphi_r = \sum_{j=1}^{n}\eta^{2j-1}_r\otimes \xi_{2j}^r - \eta^{2j}_r\otimes \xi_{2j-1}^r + \mathcal{O}_{g_0}(e^{-(a-\frac{1}{2})r})$.
Taking the limit as $r\to \infty$ shows that $\varphi = \sum_{j=1}^n \eta^{2j-1} \otimes \xi_{2j} - \eta^{2j}\otimes \xi_{2j-1}$.
Therefore, the restriction $J_0= \varphi|_{H_0}$ has at least the same regularity as $\{\eta^1|_{H_0},\ldots,\eta^{2n}|_{H_0}\}$ and $\{\xi_1,\ldots,\xi_{2n}\}$.
It follows from Theorem \ref{theorem:contact} and Lemma \ref{lemma:coframe_c1} that $J_0$ is of class $\mathcal{C}^1$.

Let us now show that $J_0$ is formally integrable.
Recall that $\gamma|_{H_0\times H_0}$ is $J_0$-invariant, so that by \cite[Proposition 5.10]{pinoy_asymptotic2021}, it suffices to show that $N_{\varphi}|_{H_0\times H_0} = \dx \eta^0|_{H_0\times H_0}\otimes \xi_0$,
where $N_{A}$ stands for the Nijenhuis tensor of the field of endomorphisms $A$, defined by $N_{A}(X,Y) = -A^2[X,Y] - [A X,AY] + A[A X,Y] + A[X,A Y]$.
Let $u$ and $v$ be any vector fields on $\dK$.
Using the fact that $\nabla$ is torsion-free, one first obtains $N_{\Phi}(Y_u,Y_v) = \Phi(\nabla_{Y_u}\Phi)Y_v - (\nabla_{\Phi Y_u}\Phi) Y_v - \Phi(\nabla_{Y_v}\Phi)Y_u + (\nabla_{\Phi Y_v}\Phi) Y_u$.
Recall that $\Phi = J - g(\cdot,\delr)\otimes \Jdelr + g(\cdot,\Jdelr)\otimes \delr$.
Since $\nabla g = 0$, $\nabla \delr = S$, $\Phi(\delr) = \Phi(\Jdelr)=0$ and $Y_u,Y_v \perp \delr$, one has
\begin{equation}
\begin{split}
\Phi(\nabla_{Y_u} \Phi)Y_v &= g(Y_v,\Jdelr)\Phi(SY_u) + \Phi(\nabla_{Y_u} J)Y_v, \\
(\nabla_{\Phi Y_u}\Phi)Y_v &= -g(Y_v,S\Phi Y_u)\Jdelr + g(Y_v,JS\Phi Y_u)\delr + g(Y_v,\Jdelr)S\Phi Y_u \\
& \quad +(\nabla_{\Phi Y_u}J)Y_v - g(Y_v,(\nabla_{\Phi Y_u}J)\delr)\delr, \\
\Phi(\nabla_{Y_v} \Phi)Y_u &= g(Y_u,\Jdelr)\Phi(SY_v) + \Phi(\nabla_{Y_v} J)Y_u, \quad \text{and}\\
(\nabla_{\Phi Y_v}\Phi)Y_u &= -g(Y_u,S\Phi Y_v)\Jdelr + g(Y_u,JS\Phi Y_v)\delr + g(Y_u,\Jdelr)S\Phi Y_v \\
& \quad + (\nabla_{\Phi Y_v}J)Y_u - g(Y_u,(\nabla_{\Phi Y_v}J)\delr)\delr.
\end{split}
\end{equation}
Recall that $\Phi$ takes values in the distribution $\{\delr\}^{\perp}$, which is involutive as the tangent field to the foliation $(\dK_r)_{r \geqslant 0}$ of $ \bar{M\setminus K}$.
The definition of the Nijenhuis tensor then shows that $N_{\Phi}$ has range in $\{\delr\}^{\perp}$.
Hence, the terms in the radial direction cancel out each others, and the remaining terms yield
\begin{equation}
\begin{split}
N_{\phi}(Y_u,Y_v) &= \left(g(Y_v,S\Phi Y_u) - g(Y_u,S\Phi Y_v)\right)\Jdelr \\
& \quad + g(Y_v,\Jdelr)\left(\Phi S Y_u - S \Phi Y_u\right) - g(Y_u,\Jdelr)\left(\Phi S Y_v - S \Phi Y_v\right) \\
& \quad + \Phi\left((\nabla_{Y_u}J)Y_v - (\nabla_{Y_v}J)Y_u\right) - \pi((\nabla_{\Phi Y_u}J)Y_v) + \pi((\nabla_{\Phi Y_v}J)Y_u)\\
&= \left(g(Y_v,S\Phi Y_u) - g(Y_u,S\Phi Y_v)\right)E_0 \\
& \quad + g(Y_v,E_0)\left(\Phi S Y_u - S \Phi Y_u\right) - g(Y_u,E_0)\left(\Phi S Y_v - S \Phi Y_v\right) \\
&\quad  + \left(g(Y_v,S\Phi Y_u) - g(Y_u,S\Phi Y_v)\right)(\Jdelr-E_0) \\
& \quad + g(Y_v,\Jdelr-E_0)\left(\Phi S Y_u - S \Phi Y_u\right) - g(Y_u,\Jdelr-E_0)\left(\Phi S Y_v - S \Phi Y_v\right) \\
& \quad + \Phi\left((\nabla_{Y_u}J)Y_v - (\nabla_{Y_v}J)Y_u\right) - \pi((\nabla_{\Phi Y_u}J)Y_v) + \pi((\nabla_{\Phi Y_v}J)Y_u),
\end{split}
\end{equation}
where $\pi$ is the orthogonal projection onto $\{\delr\}^{\perp}$.
From now, and until the rest of the proof, we assume that $u$ and $v$ are tangent to $H_0$.
Let $r \geqslant 0$, and note that $N_{\varphi_r} = \mathcal{E}_r^* N_{\Phi}$.
The \AK condition,
the uniform bound on $\|S\|_g$ (Lemma \ref{lemma:S_bounded}),
estimates on $E_0-\Jdelr$ (Corollary \ref{corollary:frame_estimates}),
estimates on $Y_u$ and $Y_v$ (Corollary \ref{corollary:Y_v_asymptotic}),
comparison between $g_0$ and $g_r$ (Corollary \ref{corollary:pinching_g_g_r}),
and estimates on $\varphi_r S_r - S_r \varphi_r$ (Lemma \ref{lemma:phi_r_estimates}),
now yield the existence of $\alpha_1 > 0$, depending on $a$, such that $N_{\varphi_r}(u,v) = e^{-r}(g(Y_v,S\Phi Y_u) - g(Y_u,S\Phi Y_v))\xi_0^r + \mathcal{O}_{g_0}(\|u\|_{g_0}\|v\|_{g_0}e^{-\alpha_1 r})$.
Similar calculations that the ones conducted to derive an expression for $\dx\eta^0_r(u,\varphi_rv)$ (see the proof of Theorem \ref{theorem:contact}) show that there exists $\alpha_2 > 0$ depending on $a$ with
\begin{equation}
e^{-r}\left(g(Y_v,S\Phi Y_u) - g(Y_u,S\Phi Y_v)\right) = \dx\eta^0(u,v) + \mathcal{O}(\|u\|_{g_0}\|v\|_{g_0}e^{-\alpha_2 r}).
\end{equation}
The $\mathcal{C}^1$ convergence of $(\varphi_r|_{H_0})_{r \geqslant 0}$ to $\varphi|_{H_0}$, and the $\mathcal{C}^0$ convergence of $(\xi_0^r)_{r \geqslant 0}$ to $\xi_0$ finally imply that $N_{\varphi}|_{H_0 \times H_0} = \lim_{r\to \infty} N_{\varphi_r}|_{H_0 \times H_0} = \dx \eta^0|_{H_0\times H_0} \otimes \xi_0$.
Consequently, $J_0$ is formally integrable.
The associated Levi-form $\dx\eta^0|_{H_0\times H_0}(\cdot,J_0\cdot)$ coincides with $\gamma|_{H_0\times H_0}$, and is thus positive definite.
Ultimately, $(\dK,H_0,J_0)$ is a strictly pseudoconvex CR manifold, which concludes the proof.
\end{proof}

\begin{remark}
If $M$ has dimension $4$, then $J_0$ is an almost complex structure of class $\mathcal{C}^1$ defined on a $2$-dimensional vector bundle.
Its integrability is automatic in this specific case.
\end{remark}

\begin{remark}
\label{remark:a>3/2_bis}
Similarly to Remark \ref{remark:a>3/2}, under the stronger assumption $a > \frac{3}{2}$, one shows that $\varphi$ is of class $\mathcal{C}^1$ in all directions.
\end{remark}


\section{The compactification}
\label{section:compactification}

\noindent We conclude this paper by showing our main Theorem.

\begin{proof}[Proof of the main Theorem]
We first give a construction for $\bar{M}$.
Fix $K$ an essential subset and $E$ its normal exponential map.
Let $M(\infty)$ be the visual boundary of $(M,g)$, which is the set of equivalent classes $[\sigma]$ of untrapped unit speed geodesic rays $\sigma$, where two rays $\sigma_1$ and $\sigma_2$ are equivalent if and only if the function $t\geqslant 0 \mapsto d_g(\sigma_1(t),\sigma_2(t))$ is bounded.
By \cite[Propositions 4.1 \& 4.4]{bahuaud_holder_2008}, $\dK$ is in bijection with $M(\infty)$ by the map $p \mapsto [E(\cdot,p)]$.
Define $\bar{M} = M \cup M(\infty)$.
The following map
\begin{equation}
\begin{array}{rccc}
\bar{\mathcal{E}}\colon &[0,1) \times \dK & \longrightarrow & \bar{M}\setminus K \\
& (\rho, p) & \longmapsto &
\begin{cases}
\mathcal{E}(-\ln \rho, p) \in M\setminus K & \text{if} \quad \rho > 0,\\
[\mathcal{E}(\cdot,p)] \in M(\infty) & \text{if} \quad \rho = 0,
\end{cases}
\end{array}
\end{equation}
is thus a bijection.
We endow $\bar{M}$ with the structure of a compact manifold with boundary through this latter bijection.
This identifies $M$ with the interior of $\bar{M}$.
Note that if $\rho > 0$, then $r = -\ln \rho$ is the distance to $K$ for $g$ in $M$.
A compactly supported modification of $\rho$ in a neighbourhood of $K$ in $M$ provides a smooth defining function for the boundary $\partial \bar{M} = M(\infty)$.
By abuse of notation, we still denote it $\rho$.

Let $\eta^0$ be the contact form and $\gamma$ be the Carnot metric given by Theorem \ref{theorem:contact}.
Let $H_0$ be the associated contact distribution, and let $J_0$ be the integrable almost complex structure on $H_0$ given by Theorem \ref{theorem:strictly_pseudoconvex}.
We see these objects as defined on $\partial\bar{M}$ through the diffeomorphism $\bar{E}(0,\cdot) \colon \{0\}\times \dK \to \partial \bar{M}$.
Then $(\partial\bar{M},H_0,J_0)$ is a strictly pseudoconvex CR manifold of class $\mathcal{C}^1$ by Theorem \ref{theorem:strictly_pseudoconvex}.
Theorem \ref{theorem:metric_expansion} and Remark \ref{remark:g-ghat} show that the metric $g$ has the desired asymptotic expansion \eqref{eq:CH_metric} near the boundary $\partial\bar{M} = \rho^{-1}(\{0\})$.

Let us show that $H_0$ and $J_0$ are induced by a continuous ambient almost complex structure $\bar{J}$.
To that end, we show that $J$ extends continuously to the boundary.
Let $\{E_0,\ldots,E_{2n}\}$ be a $J$-admissible frame on a cone $E(\R_+\times U)$, and consider the frame $ \{-\partial_{\rho}, \bar{\xi}_0,\ldots,\bar{\xi}_{2n}\}$ on $\bar{E}((0,1)\times U)$ defined by $\bar{\xi}_0 = \bar{E}^*(\rho^{-1}E_0)$ and $\bar{\xi}_j = \bar{E}^*(\rho^{-\frac{1}{2}}E_j)$ for $j\in \{1,\ldots,2n\}$.
Notice that $-\partial_{\rho} = e^{r}\delr$ on $M\setminus K$.
Proposition \ref{proposition:xi_0_estimates} and Remark \ref{remark:xi_j} show that $\{\bar{\xi}_0,\ldots,\bar{\xi}_{2n}\}$ extends continuously on the boundary $\bar{E}(\{0\}\times U)$, with limit $\{\xi_0,\ldots,\xi_{2n}\}$.
The tangent bundle of $\bar{M}$ at the boundary splits as $T\bar{M}|_{\partial\bar{M}} = \R\partial_{\rho} \oplus T\partial\bar{M} =\R\partial_{\rho} \oplus \R \xi_0 \oplus H_0$.
From the very definition of a $J$-admissible frame, one has
\begin{equation}
\begin{split}
J(e^r \delr) - e^r E_0,\quad  J(e^r E_0) + e^r \delr &= \mathcal{O}_g(e^{-(a-1)r}), \\
J(e^{\frac{r}{2}}E_{2j-1}) - e^{\frac{r}{2}}E_{2j},\quad  J(e^{\frac{r}{2}}E_{2j}) + e^{\frac{r}{2}}E_{2j-1} &= \mathcal{O}_g(e^{-(a-\frac{1}{2})r}), \quad j\in \{1,\ldots, n\}.
\end{split}
\end{equation}
It follows that in the continuous frame
$\{-\partial_{\rho},\bar{\xi}_0,\ldots,\bar{\xi}_{2n}\}$,
the matrix of $J$ reads
\begin{equation}
\left(\begin{array}{@{}c|c@{}}
\begin{matrix}
0 & -1 \\
1 & \phantom{-} 0 \\
\end{matrix}
& 0  \\
\hline
0 &
\begin{matrix}
\ddots\\
& 0 & -1 \\
& 1 & \phantom{-} 0
\end{matrix}
\end{array}\right)
+
\left(\begin{array}{@{}c|c@{}}
\begin{matrix}\\
\mathcal{O}\left(\rho^a\right) \\
\end{matrix}
& \mathcal{O}\left(\rho^{a+\frac{1}{2}}\right)  \\ \\
\hline
\mathcal{O}\left(\rho^{a-\frac{1}{2}}\right) &
\begin{matrix}
& & \phantom{a} & \\
& \mathcal{O}\left(\rho^{a} \right) & \\
& & &
\end{matrix}
\end{array}\right)
,
\end{equation}
where the top left block is of size $2\times 2$ and the bottom right block is of size $2n \times 2n$.
Hence, $J$ extends uniquely as a continuous almost complex structure $\bar{J}$ up to boundary.
In addition, $\bar{J}$ satisfies
\begin{equation}
\bar{J}(-\partial_{\rho}) = \xi_0,\quad \bar{J}\xi_0 = \partial_{\rho},\quad \bar{J}\xi_{2j-1} = \xi_{2j},\quad \text{and} \quad \bar{J}\xi_{2j} = -\xi_{2j-1},\quad j\in \{1,\ldots,2n\}.
\end{equation}
It follows that $\bar{J}|_{H_0} = J_0$, and that $H_0 = (T\partial\bar{M})\cap(\bar{J}T\partial\bar{M})$.
This concludes the proof.
\end{proof}

\begin{remark}
Under the stronger assumption that $a > \frac{3}{2}$, one can show that $\bar{J}$ is of class $\mathcal{C}^1$ up to the boundary in all directions (see Remark \ref{remark:a>3/2}).
\end{remark}

\begin{remark}
When $(M,g,J)$ is Kähler, (that is, if $\nabla J = 0$), then $(\bar{M},\bar{J})$ is a compact complex manifold with strictly pseudoconvex CR boundary.
\end{remark}



\printbibliography

\end{document}